# GLOBAL STABILIZATION OF NONLINEAR SYSTEMS BASED ON VECTOR CONTROL LYAPUNOV FUNCTIONS

**Iasson Karafyllis*** **and Zhong-Ping Jiang****

*Dept. of Environmental Eng., Technical University of Crete,
73100, Chania, Greece, email: ikarafyl@enveng.tuc.gr

**Dept. of Electrical and Computer Eng., Polytechnic Institute of New York University,
Six Metrotech Center, Brooklyn, NY 11201, U.S.A., email: zjiang@poly.edu

**Abstract**

This paper studies the use of vector Lyapunov functions for the design of globally stabilizing feedback laws for nonlinear systems. Recent results on vector Lyapunov functions are utilized. The main result of the paper shows that the existence of a vector control Lyapunov function is a necessary and sufficient condition for the existence of a smooth globally stabilizing feedback. Applications to nonlinear systems are provided: simple and easily checkable sufficient conditions are proposed to guarantee the existence of a smooth globally stabilizing feedback law. The obtained results are applied to the problem of the stabilization of an equilibrium point of a reaction network taking place in a continuous stirred tank reactor.

## 1. Introduction

Vector Lyapunov functions have been used for a long time in stability theory for nonlinear systems (see [3,6,7,11,12,13,14,15,18,19,20,21,27] and the references in [7]). Vector Lyapunov functions were first introduced by Bellman in [3] and have been acknowledged to be a more flexible tool for proving stability than the usual single Lyapunov function. The recent works [11,12,13] have established novel vector Lyapunov theorems which have qualitatively different characteristics than the classical vector Lyapunov results in [3,7, 14,15,18]. More specifically, the proof of the vector Lyapunov results in [11,12,13] is based on a small-gain analysis (rather than the classical differential inequality approach) and do not require all differential inequalities to hold simultaneously at every point of the state space.

The use of vector Lyapunov functions in control theory is not frequent. Exceptions are the works [19,20,21]. However, it seems reasonable to think that the flexibility shown by vector Lyapunov functions in stability theory can be utilized to our advantage for feedback control design in complex systems. The purpose of the present work is to show that this is indeed the case.

The main result of the present work (Theorem 3.4) is a direct extension of the well-known Artstein-Sontag theorem (see [2,5,24,25]) to the case of vector Lyapunov functions. Therefore, the term Vector Control Lyapunov Function (VCLF) is appropriate. Theorem 3.4 extends the Artstein-Sontag theorem so that

i) vector Lyapunov functions are used (instead of a single Lyapunov function),

ii) not all differential inequalities must hold at every point of the state space.

Therefore, Theorem 3.4 can allow large flexibility for its application. This feature can be crucial for feedback design in large scale systems. Recently, large scale systems have been studied intensely (see [6,8,10,22]). Section 4 of the present work shows that VCLFs can be used efficiently for the stabilization of large scale systems. Particularly, Corollaries 4.1 and 4.2 provide simple algebraic criteria that allow us to guarantee that a nonlinear system can be globally asymptotically stabilized by means of smooth feedback.

The results of the paper are applied to the global stabilization of equilibrium points of (bio)chemical reaction networks taking place in continuous stirred tank reactors (chemostats). Reaction networks have been studied in the past (see the references in [26]) and recent results have provided new insights for their properties (see [1,4,26]). Theorem 5.2 provides sufficient conditions for the existence of a smooth stabilizing feedback law for the case that the dilution rate is considered as the control input (the most frequent case in the literature). Our main assumptions on the reaction network hold for biological networks as well. For example, all chemostat models (see [23]) satisfy hypotheses (R1), (R2) and (R3) in Section 5 of the paper.

The structure of the present work is as follows. Section 2 provides the background on vector Lyapunov functions and the recent results in [12,13]. Section 3 of the present work contains the definition of the VCLF and the main result of the present work (Theorem 3.4). Section 4 is devoted to the derivation of simple sufficient conditions that guarantee the existence of smooth global stabilizers for nonlinear systems (Corollaries 4.1 and 4.2). Finally, the obtained results are applied to reaction networks in Section 5. The Appendix contains the proofs of certain auxiliary lemmas needed for the proof of Theorem 3.4.

**Notation and Terminology.** Throughout this paper we adopt the following notation and terminology:

* Let $I \subseteq \Re_+ := [0,+\infty)$ be an interval. By $L^\infty(I;U)$ (resp. $L^\infty_{loc}(I;U)$) we denote the space of measurable and (resp. locally) essentially bounded functions $u(\cdot)$ defined on $I$ and taking values in $U \subseteq \Re^m$. For a set $A \subseteq \Re^n$, $\mathrm{int}(A)$ denotes the interior of $A \subseteq \Re^n$.

* We say that a non-decreasing continuous function $\gamma : \Re_+ \to \Re_+$ is of class $N_1$ if $\gamma(0) = 0$. We say that function $\gamma : \Re_+ \to \Re_+$ is positive definite if $\gamma(0) = 0$ and $\gamma(s) > 0$ for all $s > 0$. We say that an increasing continuous function $\gamma : \Re_+ \to \Re_+$ is of class $K$ if $\gamma(0) = 0$. We say that an increasing continuous function $\gamma : \Re_+ \to \Re_+$ is of class $K_\infty$ if $\gamma(0) = 0$ and $\lim_{s \to +\infty} \gamma(s) = +\infty$. By $KL$ we denote the set of all continuous functions $\sigma = \sigma(s,t) : \Re^+ \times \Re^+ \to \Re^+$ with the properties: (i) for each $t \ge 0$ the mapping $\sigma(\cdot, t)$ is of class $K$ ; (ii) for each $s \ge 0$, the mapping $\sigma(s,\cdot)$ is non-increasing with $\lim_{t \to +\infty} \sigma(s,t) = 0$.

* For every positive integer $j$ and a set $A \subseteq \Re^n$, $C^j(A)$ ($C^j(A;\Omega)$) denotes the class of functions (taking values in $\Omega \subseteq \Re^m$) that have continuous derivatives of order $j$ on $A$ and $C^0(A;\Omega)$ denotes the class of continuous functions on $A$, which take values in $\Omega$. We also denote by $C^\infty(A) := \bigcap_{j \ge 0} C^j(A)$ the class of smooth functions on $A$.

* For a vector $x \in \Re^n$ we denote by $x'$ its transpose and by $|x|$ its Euclidean norm. $A' \in \Re^{n \times m}$ denotes the transpose of the matrix $A \in \Re^{m \times n}$. $\mathbf{1}_n$ denotes the vector $(1,...,1)' \in \mathrm{int}(\Re^n_+)$. For a vector $x \in \Re^n$ we define $\exp(x) = (\exp(x_1),...,\exp(x_n))' \in \mathrm{int}(\Re^n_+)$.

* For every scalar continuously differentiable function $V : \Re^n \to \Re$, $\nabla V(x)$ denotes the gradient of $V$ at $x \in \Re^n$, i.e., $\nabla V(x) = \left( \frac{\partial V}{\partial x_1}(x),...,\frac{\partial V}{\partial x_n}(x) \right)$. We say that a function $V : \Re^n \to \Re^+$ is positive definite if $V(x) > 0$ for all $x \ne 0$ and $V(0) = 0$. We say that a continuous function $V : \Re^n \to \Re^+$ is radially unbounded if the following property holds: "for every $M > 0$ the set $\{x \in \Re^n : V(x) \le M\}$ is compact or empty". For a vector field $f : \Re^n \to \Re^n$, $L_f V(x) = \nabla V(x) f(x)$ denotes the Lie derivative of $V$ along $f$.

* For a function $\phi : A \to \Re$, where $A \subseteq \Re^n$, $\mathrm{supp}(\phi)$ denotes the support of $\phi : A \to \Re$, i.e., $\mathrm{supp}(\phi) = \{x \in A : \phi(x) \ne 0\}$.

## 2. Background on Vector Lyapunov Functions

We consider systems described by Ordinary Differential Equations (ODEs) of the form:

$$\begin{aligned} &\dot{x} = f(d,x) \\ &x \in \Re^n ,\ d \in D \end{aligned} \tag{2.1}$$

where $D \subseteq \Re^l$ is a non-empty set and $f : D \times \Re^n \to \Re^n$ is a continuous mapping with $f(d,0) = 0$ for all $d \in D$ that satisfies the following hypotheses:

**(A1)** There exists a symmetric positive definite matrix $P \in \Re^{n \times n}$ such that for every bounded $S \subset \Re^n$, there exists a constant $L \geq 0$ satisfying the following inequality:

$$\begin{aligned} &(x-y)' P\big(f(d,x) - f(d,y)\big) \leq L|x-y|^2 \\ &\forall x \in S ,\ \forall y \in S ,\ \forall d \in D \end{aligned}$$

**(A2)** There exists $a \in K_\infty$, such that $|f(d,x)| \leq a(|x|)$ for all $(x,d) \in \Re^n \times D$.

**(A3)** There exist functions $h \in C^1(\Re^n;\Re)$ with $h(0) \leq 0$, $V_i \in C^1(\Re^n;\Re_+)$ ($i = 1,...,k$), $W \in C^1(\Re^n;\Re_+)$ being radially unbounded, a function $\delta \in C^0(\Re_+;(0,+\infty))$, a non-decreasing function $K \in C^0(\Re_+;\Re_+)$, $a_1, a_2 \in K_\infty$, $\gamma_{i,j} \in N_1$, $i, j = 1,...,k$ and a family of positive definite functions $\rho_i \in C^0(\Re_+;\Re_+)$ ($i = 1,...,k$) such that the following inequalities hold:

$$a_1(|x|) \leq \max_{i=1,...,k} V_i(x) \leq a_2(|x|), \text{ for all } x \in \Re^n \text{ satisfying } h(x) \leq 0 \tag{2.2}$$

$$\sup_{d \in D} \nabla h(x) f(d,x) \leq -\delta(h(x)) , \text{ for all } x \in \Re^n \text{ satisfying } h(x) \geq 0 \tag{2.3}$$

$$\sup_{d \in D} \nabla W(x) f(d,x) \leq K(h(x))W(x) , \text{ for all } x \in \Re^n \text{ satisfying } h(x) \geq 0 \tag{2.4}$$

Moreover, for every $i = 1,...,k$ and $x \in \Re^n$ with $h(x) \leq 0$, the following implication holds:

$$\text{“If } \max_{j=1,...,k} \gamma_{i,j}\big(V_j(x)\big) \leq V_i(x) \text{ then } \sup_{d \in D} \nabla V_i(x) f(d,x) \leq -\rho_i\big(V_i(x)\big)\text{”} \tag{2.5}$$

**Theorem 2.7 in [13]:** *Consider system (2.1) under hypotheses (A1-3). If the following set of small-gain conditions holds for each* $r = 2...,k$ *:*

$$\big(\gamma_{i_1,i_2} \circ \gamma_{i_2,i_3} \circ ... \circ \gamma_{i_r,i_1}\big)(s) < s ,\ \forall s > 0 \tag{2.6}$$

*for all* $i_j \in \{1,...,k\}$, $i_j \neq i_l$ *if* $j \neq l$, *then system (2.1) is Uniformly Robustly Globally Asymptotically Stable (URGAS) at the origin. That is, there exists a function* $\sigma \in KL$ *such that for every* $d \in L^\infty_{loc}(\Re_+;D)$, $x_0 \in \Re^n$ *the solution* $x(t) \in \Re^n$ *of (2.1) with initial condition* $x(0) = x_0$ *corresponding to input* $d \in L^\infty_{loc}(\Re_+;D)$ *satisfies*

$$|x(t)| \leq \sigma\big(|x_0|, t\big), \text{ for all } t \geq 0 \tag{2.7}$$

**Discussion and Explanations:**

**(a)** If implications (2.5) hold for all $i = 1,...,k$ and $x \in \Re^n$ then one simply takes $h(x) \equiv -1$ and arbitrary functions $W \in C^1(\Re^n;\Re_+)$, $\delta \in C^0(\Re_+;(0,+\infty))$, $K \in C^0(\Re_+;\Re_+)$. In this way we obtain Corollary 4.2 in [12]. Therefore the difference between Theorem 2.7 in [13] and Corollary 4.2 in [12] is that Theorem 2.7 assumes that the Lyapunov differential inequalities (2.5) hold only for a certain region of the state space (at the cost of the additional inequalities (2.3), (2.4)). However, Corollary 4.2 in [13] allows the study of time-varying systems as well.

**(b)** Notice that $W \in C^1(\Re^n;\Re_+)$ is not necessarily positive definite, i.e., we do not have to assume that $W(0)=0$.

**(c)** Theorem 2.7 in [13] is remarkably different from other vector Lyapunov results (see [3,7, 14,15,18]). The vector Lyapunov results in [3,7, 14,15,18] assume that the Lyapunov differential inequalities are valid on the whole state space while Hypothesis (A3) requires that each one of the differential inequalities (2.5) is valid only for a limited region of the state space (described by the inequalities $h(x) \le 0$ and $\max_{j=1,...,k} \gamma_{i,j}\left(V_j(x)\right) \le V_i(x)$). Moreover, the form of inequalities (2.5) is extremely simple and very similar to the differential inequality used for the single Lyapunov function. All these features of Theorem 2.7 in [13] will be exploited in the following section.

**(d)** Inequalities (2.3) and (2.4) allow a transient period during which the solution of (2.1) does not satisfy the Lyapunov differential inequalities (2.5). Please, see the discussion in [13].

## 3. Vector Robust Control Lyapunov Functions

Consider the feedback stabilization problem for the system:

$$\begin{gathered} \dot{x} = f(d,x) + g(x)u \\ x \in \Re^n \,,\; u \in U \subseteq \Re , d \in D \end{gathered} \tag{3.1}$$

where $D \subseteq \Re^l$ is compact and $f : D \times \Re^n \to \Re^n$, $g : \Re^n \to \Re^n$ are continuous mappings with $f(d,0)=0$ for all $d \in D$ that satisfy the following hypothesis:

**(H)** There exists a symmetric positive definite matrix $P \in \Re^{n\times n}$ such that for every bounded $S \subset \Re^n \times \Re$, there exists a constant $L \ge 0$ satisfying the following inequality:

$$\begin{gathered} (x-y)' P\big(f(d,x)+g(x)u - f(d,y) - g(y)u\big) \le L|x-y|^2 \\ \forall (x,u,y,u) \in S \times S \,,\; \forall d \in D \end{gathered}$$

For the control set $U \subseteq \Re$ we will consider the following three cases:

**(P1)** $U = \Re$.

**(P2)** There exists a constant $a \ge 0$ such that $U = [-a,+\infty)$.

**(P3)** There exist two constants $a,b \ge 0$ such that $U = [-a,b]$.

We next proceed to the definition of the Vector Robust Control Lyapunov Function (VRCLF) for system (3.1).

**Definition 3.1:** *Consider system (3.1) under hypotheses (H), (P1). Suppose that there exist functions* $\eta \in C^1(\Re^n;\Re)$ *with* $\eta(0) < 0$, $V_i \in C^1(\Re^n;\Re_+)$ *(*$i=1,...,k$*),* $W \in C^1(\Re^n;[1,+\infty))$ *being radially unbounded, a function* $\delta \in C^0(\Re_+;(0,+\infty))$*, a non-decreasing function* $K \in C^0(\Re_+;[1,+\infty))$, $\gamma_{i,j} \in N_1$, $i,j=1,...,k$, *with* $\gamma_{i,i}(s) \equiv 0$ *for* $i=1,...,k$*, a positive definite function* $\rho \in C^0(\Re_+;\Re_+)$ *and a constant* $\varepsilon > 0$ *such that the following properties hold:*

**(i)** *There exist functions* $a_1, a_2 \in K_\infty$ *such that the following inequality holds for all* $x \in \Re^n$*:*

$$a_1(|x|) \le \max_{i=1,...,k} V_i(x) \le a_2(|x|) \tag{3.2}$$

**(ii)** *The following implications hold for all* $x \in \Re^n$ *with* $\eta(x) \le \varepsilon$ *:*

$$\text{“If } \max_{j=1,...,k} \gamma_{i,j}\left(V_j(x)\right) \le V_i(x) \text{ and } L_g V_i(x) = 0 \text{ then } \max_{d \in D} L_f V_i(x) + \rho(V_i(x)) \le 0 \text{”} \tag{3.3}$$

$$\text{“If } \max_{s=1,...,k} \gamma_{i,s}\left(V_s(x)\right) \le V_i(x), \ \max_{s=1,...,k} \gamma_{j,s}\left(V_s(x)\right) \le V_j(x) \text{ and } L_g V_i(x) L_g V_j(x) < 0$$

$$\text{then } \max_{d \in D} L_f V_i(x) + \rho\left(V_i(x)\right) \le \frac{L_g V_i(x)}{L_g V_j(x)} \left( \max_{d \in D} L_f V_j(x) + \rho\left(V_j(x)\right) \right) \text{”} \tag{3.4}$$

**(iii)** *The following implications hold for all* $x \in \Re^n$ *satisfying* $\eta(x) \ge 0$ *:*

$$\text{“If } L_g \eta(x) = 0 \text{ then } \max_{d \in D} L_f \eta(x) + \delta(\eta(x)) \le 0 \text{”} \tag{3.5}$$

$$\text{“If } L_g W(x) = 0 \text{ then } \max_{d \in D} L_f W(x) - K(\eta(x)) W(x) \le 0 \text{”} \tag{3.6}$$

$$\text{“If } L_g \eta(x) L_g W(x) < 0 \text{ then } \max_{d \in D} L_f \eta(x) + \delta\left(\eta(x)\right) \le \frac{L_g \eta(x)}{L_g W(x)} \left( \max_{d \in D} L_f W(x) - K\left(\eta(x)\right) W(x) \right) \text{”} \tag{3.7}$$

**(iv)** *The following implications hold for all* $x \in \Re^n$ *with* $0 \le \eta(x) \le \varepsilon$ *:*

$$\text{“If } \max_{s=1,...k} \gamma_{j,s}\left(V_s(x)\right) \le V_j(x) \text{ and } L_g \eta(x) L_g V_j(x) < 0 \text{ then } \max_{d \in D} L_f \eta(x) + \delta\left(\eta(x)\right) \le \frac{L_g \eta(x)}{L_g V_j(x)} \left( \max_{d \in D} L_f V_j(x) + \rho\left(V_j(x)\right) \right) \text{”} \tag{3.8}$$

$$\text{“If } \max_{s=1,...k} \gamma_{j,s}\left(V_s(x)\right) \le V_j(x) \text{ and } L_g W(x) L_g V_j(x) < 0 \text{ then } \max_{d \in D} L_f W(x) - K\left(\eta(x)\right) W(x) \le \frac{L_g W(x)}{L_g V_j(x)} \left( \max_{d \in D} L_f V_j(x) + \rho\left(V_j(x)\right) \right) \text{”} \tag{3.9}$$

**(v)** *The following set of small-gain conditions holds for each* $r = 2...,k$ *:*

$$\left(\gamma_{i_1,i_2} \circ \gamma_{i_2,i_3} \circ ... \circ \gamma_{i_r,i_1}\right)(s) < s, \ \forall s > 0 \tag{3.10}$$

*for all* $i_j \in \{1,...,k\}$, $i_j \ne i_l$ *if* $j \ne l$.

**(vi)** *There exist an open set* $A \subseteq \Re^n$ *with* $0 \in A$ *and a locally Lipschitz function* $k \in C^\nu(A;U)$, *where* $\nu \in \{0,1,2,...\} \cup \{\infty\}$, *with* $k(0) = 0$, *such that the following implication holds:*

$$\text{“If } \max_{j=1,...,k} \gamma_{i,j}\left(V_j(x)\right) \le V_i(x) \text{ and } x \in A \text{ then } \max_{d \in D} L_f V_i(x) + L_g V_i(x) k(x) \le -\rho(V_i(x)) \text{”} \tag{3.11}$$

*Then, the family of functions* $V_i \in C^1(\Re^n;\Re_+)$ *(*$i = 1,...,k$*) is called a Vector Robust Control Lyapunov Function (VRCLF) for system (3.1) under hypotheses (H), (P1) or that system (3.1) under hypotheses (H), (P1) admits the Vector Robust Control Lyapunov Function* $V_i \in C^1(\Re^n;\Re_+)$ *(*$i = 1,...,k$*).*

**Definition 3.2:** *Consider system (3.1) under hypotheses (H), (P2). Suppose that there exist functions* $\eta \in C^1(\Re^n;\Re)$ *with* $\eta(0) < 0$, $V_i \in C^1(\Re^n;\Re_+)$ *(*$i = 1,...,k$*),* $W \in C^1(\Re^n;[1,+\infty))$ *being radially unbounded, a function* $\delta \in C^0(\Re_+;(0,+\infty))$, *a non-decreasing function* $K \in C^0(\Re_+;[1,+\infty))$, $\gamma_{i,j} \in N_1$, $i,j = 1,...,k$, *with* $\gamma_{i,i}(s) \equiv 0$ *for* $i = 1,...,k$, *a positive definite function* $\rho \in C^0(\Re_+;\Re_+)$ *and a constant* $\varepsilon > 0$ *such that properties (i)-(vi) hold. Moreover, assume that the following property holds:*

**(vii)** *The following implications hold:*

$$\text{“If } \max_{s=1,...,k} \gamma_{i,s}\left(V_s(x)\right) \le V_i(x), \ \eta(x) \le \varepsilon \text{ and } L_g V_i(x) > 0 \text{ then } \max_{d \in D} L_f V_i(x) + \rho\left(V_i(x)\right) - a L_g V_i(x) < 0 \text{”} \tag{3.12}$$

"If $\eta(x) \geq 0$, $L_g\eta(x) > 0$ then $\max_{d \in D} L_f\eta(x) + \delta(\eta(x)) - aL_g\eta(x) < 0$" (3.13)

"If $\eta(x) \geq 0$, $L_gW(x) > 0$ then $\max_{d \in D} L_fW(x) - K(\eta(x))W(x) - aL_gW(x) < 0$" (3.14)

*Then, the family of functions $V_i \in C^1(\Re^n;\Re_+)$ ($i = 1,...,k$) is called a Vector Robust Control Lyapunov Function (VRCLF) for system (3.1) under hypotheses (H), (P2) or that system (3.1) under hypotheses (H), (P2) admits the Vector Robust Control Lyapunov Function $V_i \in C^1(\Re^n;\Re_+)$ ($i = 1,...,k$).*

**Definition 3.3:** *Consider system (3.1) under hypotheses (H), (P3). Suppose that there exist functions $\eta \in C^1(\Re^n;\Re)$ with $\eta(0) < 0$, $V_i \in C^1(\Re^n;\Re_+)$ ($i = 1,...,k$), $W \in C^1(\Re^n;[1,+\infty))$ being radially unbounded, a function $\delta \in C^0(\Re_+;(0,+\infty))$, a non-decreasing function $K \in C^0(\Re_+;[1,+\infty))$, $\gamma_{i,j} \in N_1$, $i, j = 1,...,k$, with $\gamma_{i,i}(s) \equiv 0$ for $i = 1,...,k$, a positive definite function $\rho \in C^0(\Re_+;\Re_+)$ and a constant $\varepsilon > 0$ such that properties (i)-(vii) hold. Moreover, assume that the following property holds:*

**(viii)** *The following implications hold:*

"If $\max_{s=1,...,k} \gamma_{i,s}(V_s(x)) \leq V_i(x)$, $\eta(x) \leq \varepsilon$ and $L_gV_i(x) < 0$ then $\max_{d \in D} L_fV_i(x) + \rho(V_i(x)) + bL_gV_i(x) < 0$" (3.15)

"If $\eta(x) \geq 0$, $L_g\eta(x) < 0$ then $\max_{d \in D} L_f\eta(x) + \delta(\eta(x)) + bL_g\eta(x) < 0$" (3.16)

"If $\eta(x) \geq 0$, $L_gW(x) < 0$ then $\max_{d \in D} L_fW(x) - K(\eta(x))W(x) + bL_gW(x) < 0$" (3.17)

*Then, the family of functions $V_i \in C^1(\Re^n;\Re_+)$ ($i = 1,...,k$) is called a Vector Robust Control Lyapunov Function (VRCLF) for system (3.1) under hypotheses (H), (P3) or that system (3.1) under hypotheses (H), (P3) admits the Vector Robust Control Lyapunov Function $V_i \in C^1(\Re^n;\Re_+)$ ($i = 1,...,k$).*

We are now ready to state the main result of the section.

**Theorem 3.4:** *Consider system (3.1) under hypothesis (H) and under one of the hypotheses (P1), (P2), (P3). Suppose that system (3.1) admits the VRCLF $V_i \in C^1(\Re^n;\Re_+)$ ($i = 1,...,k$). Then, there exists a locally Lipschitz function $\tilde{k} \in C^\nu(\Re^n;U) \cap C^\infty(\Re^n \setminus \{0\};U)$, where $\nu \in \{0,1,2,...\} \cup \{\infty\}$ is the number involved in property (vi) of the definition of the VRCLF, with $\tilde{k}(0) = 0$, such that the closed-loop system (3.1) with $u = \tilde{k}(x)$ is URGAS.*

Here, it should be noted that the existence of a VRCLF is a necessary condition as well for the existence of a smooth globally stabilizing feedback if the vector fields $f : D \times \Re^n \to \Re^n$, $g : \Re^n \to \Re^n$ are locally Lipschitz. Indeed, the converse Lyapunov theorem in [16] guarantees the existence of a single Lyapunov function for the closed-loop system which satisfies all requirements of Definition 3.1 (or Definitions 3.2, 3.3) with $k = 1$.

For the proof of Theorem 3.4 we need two technical lemmas. The first lemma deals with a set of inequalities. Its proof is given in the Appendix.

**Lemma 3.5:** *Let $a, b \geq 0$ with $a + b > 0$, $f_i, g_i \in \Re$ ($i = 1,...,m$) be given constants. There exists $u \in \Re$ such that $f_i + g_iu < 0$, for all $i = 1,...,m$ if and only if the following implications hold:*

**(I)** If $g_i = 0$ for certain $i \in \{1,...,m\}$, then $f_i < 0$,

**(II)** If $g_ig_j < 0$ for certain pair $(i, j) \in \{1,...,m\} \times \{1,...,m\}$, then $f_i < \frac{g_i}{g_j} f_j$.

*Moreover, there exists $u \in (-a,+\infty)$ such that $f_i + g_i u < 0$, for all $i = 1,...,m$ if and only if implications (I), (II) hold and the following implication holds as well:*

**(III)** If $g_i > 0$ for certain $i \in \{1,...,m\}$, then $f_i - a g_i < 0$.

*Finally, there exists $u \in (-a,b)$ such that $f_i + g_i u < 0$, for all $i = 1,...,m$ if and only if implications (I), (II), (III) hold and the following implication holds as well:*

**(IV)** If $g_i < 0$ for certain $i \in \{1,...,m\}$, then $f_i + b g_i < 0$.

The second technical lemma guarantees that we may assume that all requirements (i)-(vi) for a VRCLF hold for functions $\gamma_{i,j} \in N_1$, $i,j = 1,...,k$, which are positive definite for $i \neq j$. Its proof is provided at the Appendix.

**Lemma 3.6:** *Suppose that system (3.1) admits the VRCLF $V_i \in C^1(\Re^n;\Re_+)$ ($i = 1,...,k$). Then all properties of the definition of the VRCLF hold with functions $\gamma_{i,j} \in N_1$, $i,j = 1,...,k$, which are positive definite for $i \neq j$ and satisfy $\lim_{s \to +\infty} \gamma_{i,j}(s) = +\infty$ for all $i,j = 1,...,k$ with $i \neq j$.*

We are now ready to prove Theorem 3.4.

**Proof of Theorem 3.4:** Without loss of generality and since $\eta(0) < 0$ we may assume that the neighborhood $A$ involved in hypothesis (vi) satisfies $A \subset \left\{ x \in \Re^n : \eta(x) < 0 \right\}$. Let $r > 0$ with $\left\{ x \in \Re^n : |x| \leq 2r \right\} \subset A$. Moreover, by virtue of Lemma 3.6, without loss of generality we may assume that all functions $\gamma_{i,j} \in N_1$, $i,j = 1,...,k$, are positive definite for $i \neq j$.

Using hypotheses (ii), (iii), (iv), convexity of $U \subseteq \Re$ and partition of unity arguments, we will next construct smooth feedback laws $k_1 : \left\{ x \in \Re^n : \eta(x) > 0 \right\} \to U$, $k_2 : \left\{ x \in \Re^n : \eta(x) < \varepsilon, x \neq 0 \right\} \to U$ and $k_3 : \left\{ x \in \Re^n : 0 < \eta(x) < \varepsilon \right\} \to U$ such that the following hold:

$$\max_{d \in D} L_f \eta(x) + k_1(x) L_g \eta(x) + \frac{1}{2} \delta(\eta(x)) \leq 0, \text{ for all } x \in \left\{ x \in \Re^n : \eta(x) > 0 \right\} \tag{3.18}$$

$$\max_{d \in D} L_f W(x) + k_1(x) L_g W(x) - 2K(\eta(x)) W(x) \leq 0, \text{ for all } x \in \left\{ x \in \Re^n : \eta(x) > 0 \right\} \tag{3.19}$$

$$\text{If } \max_{j=1,...,k} \gamma_{i,j}\left(V_j(x)\right) \leq V_i(x),\ x \neq 0 \text{ and } \eta(x) < \varepsilon \text{ then } \max_{d \in D} L_f V_i(x) + k_2(x) L_g V_i(x) + \frac{1}{2} \rho(V_i(x)) \leq 0 \tag{3.20}$$

$$\max_{d \in D} L_f \eta(x) + k_3(x) L_g \eta(x) + \frac{1}{2} \delta(\eta(x)) \leq 0, \text{ for all } x \in \left\{ x \in \Re^n : 0 < \eta(x) < \varepsilon \right\} \tag{3.21}$$

$$\max_{d \in D} L_f W(x) + k_3(x) L_g W(x) - 2K(\eta(x)) W(x) \leq 0, \text{ for all } x \in \left\{ x \in \Re^n : 0 < \eta(x) < \varepsilon \right\} \tag{3.22}$$

$$\text{If } \max_{j=1,...,k} \gamma_{i,j}\left(V_j(x)\right) \leq V_i(x),\ x \neq 0 \text{ and } 0 < \eta(x) < \varepsilon, \text{ then } \max_{d \in D} L_f V_i(x) + k_3(x) L_g V_i(x) + \frac{1}{2} \rho(V_i(x)) \leq 0 \tag{3.23}$$

Let $p : \Re \to [0,1]$ a smooth non-decreasing function with $p(s) = 0$ for all $s \leq 0$ and $p(s) = 1$ for all $s \geq 1$. We define:

$$\tilde{k}(x) := k_1(x) \text{ for all } x \in \left\{ x \in \Re^n : \eta(x) > \frac{4\varepsilon}{5} \right\} \tag{3.24}$$

$$\tilde{k}(x) := k_2(x) \text{ for all } x \in \left\{ x \in \Re^n : \eta(x) < \frac{\varepsilon}{5} \right\} \cap \left\{ x \in \Re^n : |x| > 2r \right\} \tag{3.25}$$

$$\tilde{k}(x) := k_3(x) \text{ for all } x \in \left\{ x \in \Re^n : \frac{2\varepsilon}{5} < \eta(x) < \frac{3\varepsilon}{5} \right\} \tag{3.26}$$

$$\tilde{k}(x) := \left(1 - p\left(\frac{5\eta(x)}{\varepsilon} - 3\right)\right) k_3(x) + p\left(\frac{5\eta(x)}{\varepsilon} - 3\right) k_1(x) \text{ for all } x \in \left\{ x \in \Re^n : \frac{3\varepsilon}{5} \le \eta(x) \le \frac{4\varepsilon}{5} \right\} \tag{3.27}$$

$$\tilde{k}(x) := \left(1 - p\left(\frac{5\eta(x)}{\varepsilon} - 1\right)\right) k_2(x) + p\left(\frac{5\eta(x)}{\varepsilon} - 1\right) k_3(x) \text{ for all } x \in \left\{ x \in \Re^n : \frac{\varepsilon}{5} \le \eta(x) \le \frac{2\varepsilon}{5} \right\} \tag{3.28}$$

$$\tilde{k}(x) := k(x) \text{ for all } x \in \left\{ x \in \Re^n : |x| < r \right\} \tag{3.29}$$

$$\tilde{k}(x) := \left(1 - p\left(\frac{|x|^2 - r^2}{3r^2}\right)\right) k(x) + p\left(\frac{|x|^2 - r^2}{3r^2}\right) k_2(x) \text{ for all } x \in \left\{ x \in \Re^n : r \le |x| \le 2r \right\} \tag{3.30}$$

where $k \in C^{\nu}(A;U)$ is the locally Lipschitz mapping involved in hypothesis (vi).

It follows directly from hypothesis (H), compactness of $D \subseteq \Re^l$ and the fact that the mapping $\tilde{k} : \Re^n \to U$ defined above is locally Lipschitz with $\tilde{k}(0) = 0$, that the closed-loop system (3.1) with $u = \tilde{k}(x)$ satisfies hypotheses (A1-3) with $h(x) := \eta(x) - \frac{2\varepsilon}{5}$, $\rho_i(s) := \frac{1}{2}\rho(s)$ ($i = 1,...,k$) and $\tilde{\delta}(s) := \frac{1}{2}\delta\left(s + \frac{2\varepsilon}{5}\right)$, $\tilde{K}(s) := 2K\left(s + \frac{2\varepsilon}{5}\right)$ in place of $\delta \in C^0(\Re_+;(0,+\infty))$, $K \in C^0(\Re_+;[1,+\infty))$. Consequently, we may conclude from Theorem 2.7 in [13] that the closed-loop system (3.1) with $u = \tilde{k}(x)$ is URGAS.

Therefore, we are left with the task of constructing smooth feedback laws $k_1 : \{x \in \Re^n : \eta(x) > 0\} \to U$, $k_2 : \{x \in \Re^n : \eta(x) < \varepsilon, x \neq 0\} \to U$ and $k_3 : \{x \in \Re^n : 0 < \eta(x) < \varepsilon\} \to U$ so that (3.18)-(3.23) hold.

Construction of $k_1 : \{x \in \Re^n : \eta(x) > 0\} \to U$.

By virtue of (3.5), (3.6), (3.7), (3.13), (3.14), (3.16), (3.17) and Lemma 3.5, it follows that for every $x \in \Re^n$ with $\eta(x) > 0$ there exists $u_x \in U$ with $\max_{d \in D} L_f\eta(x) + u_x L_g\eta(x) + \frac{3}{4}\delta(\eta(x)) \le 0$ and $\max_{d \in D} L_f W(x) + u_x L_g W(x) - \frac{3}{2}K(\eta(x))W(x) \le 0$. Continuity of $\max_{d \in D} L_f\eta(x), L_g\eta(x), \delta(\eta(x))$, $\max_{d \in D} L_f W(x), L_g W(x), K(\eta(x)), W(x)$ implies that there exists $\delta_x > 0$ such that $\eta(y) > 0$, $\max_{d \in D} L_f\eta(y) + u_x L_g\eta(y) + \frac{1}{2}\delta(\eta(y)) \le 0$, $\max_{d \in D} L_f W(y) + u_x L_g W(y) - 2K(\eta(y))W(y) \le 0$ for all $y \in \Re^n$ with $|y - x| \le \delta_x$.

Therefore, the sets $\{y \in \Re^n : |y - x| < \delta_x\}$ for all $x \in \Re^n$ with $\eta(x) > 0$ form an open covering of the open set $\{x \in \Re^n : \eta(x) > 0\}$. By partition of unity, there exists a family of smooth functions $\phi_i : \{x \in \Re^n : \eta(x) > 0\} \to [0,1]$, $i = 1,2,...$ such that:

a) For every $i = 1,2,...$ there exists $x_i \in \Re^n$ with $\eta(x_i) > 0$ and $\text{supp}(\phi_i) \subset \{y \in \Re^n : |y - x_i| < \delta_{x_i}\}$.

b) The sum $\sum_i \phi_i(x)$ is locally finite and satisfies $\sum_{i=1}^{\infty} \phi_i(x) = 1$ for all $x \in \Re^n$ with $\eta(x) > 0$.

We define for all $x \in \Re^n$ with $\eta(x) > 0$:

$$k_1(x) := \sum_{i=1}^{\infty} u_{x_i} \phi_i(x) \tag{3.31}$$

We first notice that (by local finiteness) the mapping $k_1$ defined by (3.31) is smooth and satisfies $k_1(x) \in U$ for all $x \in \Re^n$ with $\eta(x) > 0$.

For arbitrary $x \in \Re^n$ with $\eta(x) > 0$ we define the finite set $J_x \subset \{1,2,...\}$ of indices such that $\phi_i(x) > 0$. It follows that $x \in \left\{ y \in \Re^n : \left| y - x_j \right| < \delta_{x_j} \right\}$ for all $j \in J_x$. Therefore, we get $\max_{d \in D} L_f \eta(x) + u_{x_j} L_g \eta(x) + \frac{1}{2} \delta(\eta(x)) \le 0$, $\max_{d \in D} L_f W(x) + u_{x_j} L_g W(x) - 2K(\eta(x)) W(x) \le 0$ for all $j \in J_x$. Using the previous inequalities, the fact that $\sum_{j \in J_x} \phi_j(x) = 1$ and definition (3.31) we obtain (3.18) and (3.19).

Construction of $k_2 : \left\{ x \in \Re^n : \eta(x) < \varepsilon, x \neq 0 \right\} \to U$.

Let $x \in \Re^n$ with $\eta(x) < \varepsilon$, $x \neq 0$ be arbitrary. Le $J^+(x)$ be the set of all $j \in \{1,...,k\}$ such that $\max_{s=1,...,k} \gamma_{j,s}(V_s(x)) \le V_j(x)$. Since all $\gamma_{i,j} \in N_1$, $i, j = 1,...,k$ are positive definite, we have $V_j(x) > 0$ for all $j \in J^+(x)$. By virtue of (3.3), (3.4), (3.12), (3.15) and Lemma 3.5, there exists $u_x \in U$ with

$$\max_{d \in D} L_f V_j(x) + u_x L_g V_j(x) + \frac{3}{4} \rho(V_j(x)) \le 0 \text{, for all } j \in J^+(x) \tag{3.32}$$

Define $J^-(x)$ to be the set of all $j \in \{1,...,k\}$ such that $\max_{s=1,...,k} \gamma_{j,s}(V_s(x)) > V_j(x)$. Notice that $J^+(x) \cup J^-(x) = \{1,...,k\}$ for all $x \in \Re^n$ with $\eta(x) < \varepsilon$, $x \neq 0$.

Continuity of mappings $\max_{d \in D} L_f V_j(x), L_g V_j(x), \rho(V_j(x))$ ($j = 1,...,k$), $\gamma_{i,j} \in N_1$, $i, j = 1,...,k$ implies that there exists $\delta_x > 0$ such that $\eta(y) < \varepsilon$, $y \neq 0$, $\max_{s=1,...,k} \gamma_{j,s}(V_s(y)) > V_j(y)$ for all $j \in J^-(x)$, $\max_{d \in D} L_f V_j(y) + u_x L_g V_j(y) + \frac{1}{2} \rho(V_j(y)) \le 0$ for all $j \in J^+(x)$ and for all $y \in \Re^n$ with $|y - x| \le \delta_x$. Notice that the previous inequality imply that $J^+(y) \subseteq J^+(x)$ for all $y \in \Re^n$ with $|y - x| \le \delta_x$. Therefore, we have

$$\max_{d \in D} L_f V_j(y) + u_x L_g V_j(y) + \frac{1}{2} \rho(V_j(y)) \le 0 \text{, for all } j \in J^+(y) \text{ and } y \in \Re^n \text{ with } |y - x| \le \delta_x \tag{3.33}$$

Therefore, the sets $\left\{ y \in \Re^n : |y - x| < \delta_x \right\}$ for all $x \in \Re^n$ with $\eta(x) < \varepsilon$, $x \neq 0$ form an open covering of the open set $\left\{ x \in \Re^n : \eta(x) < \varepsilon, x \neq 0 \right\}$. By partition of unity, there exists a family of smooth functions $\phi_i : \left\{ x \in \Re^n : \eta(x) < \varepsilon, x \neq 0 \right\} \to [0,1]$, $i = 1,2,...$ such that:

a) For every $i = 1,2,...$ there exists $x_i \in \Re^n$ with $\eta(x_i) < \varepsilon$, $x_i \neq 0$ and $\mathrm{supp}(\phi_i) \subset \left\{ y \in \Re^n : \left| y - x_i \right| < \delta_{x_i} \right\}$.

b) The sum $\sum_i \phi_i(x)$ is locally finite and satisfies $\sum_{i=1}^{\infty} \phi_i(x) = 1$ for all $x \in \Re^n$ with $\eta(x) < \varepsilon$, $x \neq 0$.

We define for all $x \in \Re^n$ with $\eta(x) < \varepsilon$, $x \neq 0$:

$$k_2(x) := \sum_{i=1}^{\infty} u_{x_i} \phi_i(x) \tag{3.34}$$

We first notice that (by local finiteness) the mapping $k_2$ defined by (3.34) is smooth and satisfies $k_2(x) \in U$ for all $x \in \Re^n$ with $\eta(x) < \varepsilon$, $x \neq 0$.

For arbitrary $x \in \Re^n$ with $\eta(x) < \varepsilon$, $x \neq 0$ we define the finite set $I_x \subset \{1,2,...\}$ of indices such that $\phi_i(x) > 0$. It follows that $x \in \left\{ y \in \Re^n : |y - x_i| < \delta_{x_i} \right\}$ for all $i \in I_x$. Therefore, by virtue of (3.33), the fact that $\sum_{i \in I_x} \phi_i(x) = 1$, the fact that $J^+(x)$ is the set of all $j \in \{1,...,k\}$ such that $\max_{s=1,...,k} \gamma_{j,s}(V_s(x)) \leq V_j(x)$ and definition (3.34), we obtain (3.20).

The construction of $k_3 : \left\{ x \in \Re^n : 0 < \eta(x) < \varepsilon \right\} \to U$ is similar to the constructions of $k_1 : \left\{ x \in \Re^n : \eta(x) > 0 \right\} \to U$ and $k_2 : \left\{ x \in \Re^n : \eta(x) < \varepsilon, x \neq 0 \right\} \to U$ (using partition of unity arguments). The proof is complete. $\lhd$

It should be noted that a different way of proving Theorem 3.4 is by using Michael's theorem (instead of partition of unity arguments). However, the use of Michael's theorem results to a continuous feedback (instead of a smooth feedback law).

## 4. Application to Feedback Stabilizer Design

A natural question of practical importance is how to construct a VRCLF which satisfies the (involved) assumptions of the VRCLF. The purpose of this section is to provide conditions that allow us to use simple VRCLFs (namely, the functions $V_i(x) = \frac{1}{2} x_i^2$ for $i = 1,...,n$). More specifically, this section is devoted to the feedback stabilization problem for nonlinear systems of the form:

$$\dot{x}_i = f_i(x) + g_i(x)u \,, \; i = 1,..,n \tag{4.1}$$

where $x = (x_1,...,x_n)' \in \Re^n$, $u \in U \subseteq \Re$ and the mappings $f_i : \Re^n \to \Re$, $g_i : \Re^n \to \Re$ are locally Lipschitz with $f_i(0) = 0$ for $i = 1,..,n$. More specifically, we will use the results of the previous section in order to derive sufficient conditions for the existence of a globally stabilizing feedback. Interestingly, the developed sufficient conditions will be easily checkable even for large $n$ (large scale systems).

Our main results are the following corollaries. Their proof is based on Theorem 3.4.

**Corollary 4.1:** *Consider system (4.1) under hypothesis (P1) and suppose that there exist functions $\eta \in C^1(\Re^n;\Re)$ with $\eta(0) < 0$, $W \in C^1(\Re^n;[1,+\infty))$ being radially unbounded, a function $\delta \in C^0(\Re_+;(0,+\infty))$, a non-decreasing function $K \in C^0(\Re_+;[1,+\infty))$, $\tilde{\gamma}_{i,j} \in N_1$, $i,j = 1,...,n$, with $\gamma_{i,i}(s) \equiv 0$ for $i = 1,...,n$ which satisfy the small-gain conditions (3.10), a function $Q \in C^0(\Re;\Re_+)$ with $Q(x) > 0$ for all $x \neq 0$ and a constant $\varepsilon > 0$ such that the following implications hold for all $i,j \in \{1,...,n\}$:*

$$\text{“If } \max_{s=1,...,n} \tilde{\gamma}_{i,s}(|x_s|) \leq |x_i| \,,\; \max_{s=1,...,n} \tilde{\gamma}_{j,s}(|x_s|) \leq |x_j| \,,\; \eta(x) \leq \varepsilon \text{ and } x_i x_j g_i(x) g_j(x) < 0$$

$$\text{then } \frac{f_i(x) g_j(x) - f_j(x) g_i(x)}{x_i Q(x_i) g_j(x) - x_j Q(x_j) g_i(x)} \leq -1 \text{”} \tag{4.2}$$

$$\text{“If } g_i(x) = 0 \,,\; \eta(x) \leq \varepsilon \text{ and } \max_{s=1,...,n} \tilde{\gamma}_{i,s}(|x_s|) \leq |x_i| \text{ then } x_i f_i(x_i) + x_i^2 Q(x_i) \leq 0 \text{”} \tag{4.3}$$

$$\text{“If } \max_{s=1,...,n} \tilde{\gamma}_{j,s}(|x_s|) \leq |x_j| \,,\; 0 \leq \eta(x) \leq \varepsilon \text{ and } x_j g_j(x) \sum_{l=1}^{n} \frac{\partial \eta}{\partial x_l}(x) g_l(x) < 0$$

$$\text{then } \sum_{l=1}^{n} \frac{\partial \eta}{\partial x_l}(x) \left( f_l(x) - g_l(x) \frac{f_j(x) + x_j Q(x_j)}{g_j(x)} \right) \leq -\delta(\eta(x)) \text{”} \tag{4.4}$$

"If $\max_{s=1,...,n} \tilde{\gamma}_{j,s}(|x_s|) \le |x_j|$, $0 \le \eta(x) \le \varepsilon$ and $x_j g_j(x) \sum_{l=1}^{n} \frac{\partial W}{\partial x_l}(x) g_l(x) < 0$

then $\sum_{l=1}^{n} \frac{\partial W}{\partial x_l}(x) \left( f_l(x) - g_l(x) \frac{f_j(x) + x_j Q(x_j)}{g_j(x)} \right) \le K(\eta(x)) W(x)$" (4.5)

*Moreover, suppose that there exists a vector* $k \in \Re^n$ *such that* $x_i f_i(x) + x_i g_i(x) k'x \le -x_i^2 Q(x_i)$ *for all* $x = (x_1,...,x_n)' \in \Re^n$ *in a neighborhood of* $0 \in \Re^n$ *with* $\max_{s=1,...,k} \tilde{\gamma}_{i,s}(|x_s|) \le |x_i|$ *(* $i = 1,...,n$ *). Finally, suppose that there exists a function* $\bar{k} \in C^0(\Re^n;\Re)$ *such that* $\sum_{i=1}^{n} \frac{\partial W}{\partial x_i}(x) \left( f_i(x) + g_i(x)\bar{k}(x) \right) \le K(\eta(x)) W(x)$ *and* $\sum_{i=1}^{n} \frac{\partial \eta}{\partial x_i}(x) \left( f_i(x) + g_i(x)\bar{k}(x) \right) \le -\delta(\eta(x))$ *for all* $x \in \Re^n$ *with* $\eta(x) \ge 0$.

*Then there exists* $\tilde{k} \in C^\infty(\Re^n;\Re)$, *with* $\tilde{k}(0) = 0$, *such that* $0 \in \Re^n$ *is UGAS for the closed-loop system* $\dot{x}_i = f_i(x) + g_i(x)\tilde{k}(x)$, $i = 1,..,n$.

**Corollary 4.2:** *Consider system (4.1) under hypothesis (P3) with* $a,b > 0$ *and suppose that there exist functions* $\eta \in C^1(\Re^n;\Re)$ *with* $\eta(0) < 0$, $W \in C^1(\Re^n;[1,+\infty))$ *being radially unbounded, a function* $\delta \in C^0(\Re_+;(0,+\infty))$, *a non-decreasing function* $K \in C^0(\Re_+;[1,+\infty))$, $\tilde{\gamma}_{i,j} \in N_1$, $i,j = 1,...,n$, *with* $\gamma_{i,i}(s) \equiv 0$ *for* $i = 1,...,n$ *which satisfy the small-gain conditions (3.10), a function* $Q \in C^0(\Re;\Re_+)$ *with* $Q(x) > 0$ *for all* $x \ne 0$ *and a constant* $\varepsilon > 0$ *such that implications (4.2), (4.3), (4.4) and (4.5) hold for all* $i,j \in \{1,...,n\}$. *Moreover, suppose that there exists a vector* $k \in \Re^n$ *such that* $x_i f_i(x) + x_i g_i(x) k'x \le -x_i^2 Q(x_i)$ *for all* $x = (x_1,...,x_n)' \in \Re^n$ *in a neighborhood of* $0 \in \Re^n$ *with* $\max_{s=1,...,k} \tilde{\gamma}_{i,s}(|x_s|) \le |x_i|$ *(* $i = 1,...,n$ *) and suppose that there exists a function* $\bar{k} \in C^0(\Re^n;U)$ *such that* $\sum_{i=1}^{n} \frac{\partial W}{\partial x_i}(x) \left( f_i(x) + g_i(x)\bar{k}(x) \right) \le K(\eta(x)) W(x)$ *and* $\sum_{i=1}^{n} \frac{\partial \eta}{\partial x_i}(x) \left( f_i(x) + g_i(x)\bar{k}(x) \right) \le -\delta(\eta(x))$ *for all* $x \in \Re^n$ *with* $\eta(x) \ge 0$. *Finally, suppose that the following implication holds:*

"If $\max_{s=1,...,n} \tilde{\gamma}_{i,s}(|x_s|) \le |x_i|$, $\eta(x) \le \varepsilon$ and $x_i g_i(x) > 0$ then $x_i f_i(x) + x_i^2 Q(x_i) - a x_i g_i(x) < 0$" (4.6)

"If $\max_{s=1,...,n} \tilde{\gamma}_{i,s}(|x_s|) \le |x_i|$, $\eta(x) \le \varepsilon$ and $x_i g_i(x) < 0$ then $x_i f_i(x) + x_i^2 Q(x_i) + b x_i g_i(x) < 0$" (4.7)

*Then there exists* $\tilde{k} \in C^\infty(\Re^n;U)$, *with* $\tilde{k}(0) = 0$, *such that* $0 \in \Re^n$ *is UGAS for the closed-loop system* $\dot{x}_i = f_i(x) + g_i(x)\tilde{k}(x)$, $i = 1,..,n$.

**Proof of Corollary 4.1 and Corollary 4.2:** A direct application of Theorem 3.4 with $k = n$, $V_i(x) = \frac{1}{2}x_i^2$ for $i = 1,...,n$, $\gamma_{i,j}(s) := \frac{1}{2}\left(\tilde{\gamma}_{i,j}\left(\sqrt{2s}\right)\right)^2$ for $i,j = 1,...,n$ and $\rho(s) := 2s \min\left(Q\left(\sqrt{2s}\right), Q\left(-\sqrt{2s}\right)\right)$, $a_1(s) := \frac{1}{2n}s^2$, $a_2(s) := \frac{1}{2}s^2$ for all $s \ge 0$.

Notice that implications (3.3), (3.4) are directly implied by implications (4.3), (4.2), respectively and the above definitions. Implications (3.5), (3.6), (3.7) are a direct consequence of the existence of a function $\bar{k} \in C^0(\Re^n;U)$ such that $\sum_{i=1}^{n}\frac{\partial W}{\partial x_i}(x)\left(f_i(x)+g_i(x)\bar{k}(x)\right)\le K\left(\eta(x)\right)W(x)$ and $\sum_{i=1}^{n}\frac{\partial \eta}{\partial x_i}(x)\left(f_i(x)+g_i(x)\bar{k}(x)\right)\le -\delta\left(\eta(x)\right)$ for all $x\in\Re^n$ with $\eta(x)\ge 0$.

For Corollary 4.2, we notice that (3.12), (3.15) are equivalent to implications (4.6), (4.7) and that implications (3.13), (3.14), (3.16), (3.17) are directly implied by the inequalities $\sum_{i=1}^{n}\frac{\partial W}{\partial x_i}(x)\left(f_i(x)+g_i(x)\bar{k}(x)\right)\le K\left(\eta(x)\right)W(x)$ and $\sum_{i=1}^{n}\frac{\partial \eta}{\partial x_i}(x)\left(f_i(x)+g_i(x)\bar{k}(x)\right)\le -\delta\left(\eta(x)\right)$ for all $x\in\Re^n$ with $\eta(x)\ge 0$ (possibly by replacing $K\left(\eta(x)\right)$ and $\delta\left(\eta(x)\right)$ by $2K\left(\eta(x)\right)$ and $\frac{1}{2}\delta\left(\eta(x)\right)$, respectively) and the fact that $\bar{k}(x)\in U$ for all $x\in\Re^n$. Notice that since $a,b>0$, hypothesis (vi) of Definition 3.2 holds with the linear feedback law $u=k'x$ (possibly by replacing the initial neighborhood by another neighborhood with $|k||x|\le\min(a,b)$). The proof is complete. ◁

A direct application of Corollary 4.1 and Corollary 4.2 is obtained by selecting $\eta(x)\equiv -1$, $\varepsilon=-1$: in this case implications (4.4), (4.5) are automatically satisfied for arbitrary functions $W\in C^1(\Re^n;[1,+\infty))$, $\delta\in C^0(\Re_+;(0,+\infty))$, $K\in C^0(\Re_+;[1,+\infty))$. Moreover, the assumption of the existence of a function $\bar{k}\in C^0\left(\Re^n;U\right)$ such that $\sum_{i=1}^{n}\frac{\partial W}{\partial x_i}(x)\left(f_i(x)+g_i(x)\bar{k}(x)\right)\le K\left(\eta(x)\right)W(x)$ and $\sum_{i=1}^{n}\frac{\partial \eta}{\partial x_i}(x)\left(f_i(x)+g_i(x)\bar{k}(x)\right)\le -\delta\left(\eta(x)\right)$ for all $x\in\Re^n$ with $\eta(x)\ge 0$ is automatically satisfied for arbitrary $\bar{k}\in C^0\left(\Re^n;U\right)$. The reader should notice that implications (4.2), (4.3) are easily checkable: the gain functions $\tilde{\gamma}_{i,j}\in N_1$, $i,j=1,...,n$ are selected so that implications (4.2), (4.3) hold. The following example shows how to use Corollary 4.1 with $\eta(x)\equiv -1$, $\varepsilon=-1$ for a third-order nonlinear system.

**Example 4.3:** Consider the nonlinear system

$$\begin{aligned}\dot{x}_1 &= -x_1+x_2\\ \dot{x}_2 &= -x_2+g(x)u\\ \dot{x}_3 &= x_1^2+u\\ x &= (x_1,x_2,x_3)'\in\Re^3, u\in\Re\end{aligned} \tag{4.8}$$

where $g:\Re^3\to\Re$ is a locally Lipschitz function. The problem that we study in this example is:

**(Q)** "For what functions $g:\Re^3\to\Re$, system (4.8) can be stabilized globally by a smooth feedback?"

Clearly, (4.8) is stabilizable by means of a smooth feedback if $g(x)\equiv 0$. Using the function $V(x)=\frac{1}{2}x_1^2+\frac{p}{2}x_2^2+\frac{q}{2}x_3^2$ as a CLF candidate, where $p,q>0$, we conclude that this function is a CLF for system (4.8) provided that the following condition holds:

$$\text{For each } x\in\Re^3\setminus\{0\} \text{ satisfying } qx_3=-px_2g(x) \text{ it holds that } x_1^2-x_1x_2+px_2^2>qx_1^2x_3 \tag{4.9}$$

Here, we obtain different conditions for the function $g:\Re^3\to\Re$, which allow smooth stabilizability of system (4.8). We show that system (4.8) can be stabilized globally by smooth feedback provided that there exist constants $\lambda\in(0,1)$, $\sigma\in(0,1)$, functions $\gamma\in K_\infty\cap C^1(\Re_+)$ with $\frac{d\gamma}{ds}(0)>\frac{\sigma|g(0)|}{1-\sigma}$, $Q\in C^0(\Re;\Re_+)$ with $0<Q(x)\le\sigma$ for all $x\neq 0$, such that the locally Lipschitz function $g:\Re^3\to\Re$ satisfies the following implication:

$$\lambda(1-\sigma)|x_1|\le\gamma\left(|x_3|\right),\ \lambda|x_2|\le\gamma\left(|x_3|\right)\le|x_2|,\ \ x_2x_3g(x)<0\Rightarrow x_2^2+x_2x_1^2g(x)\ge|x_2x_3g(x)|Q(x_3)+x_2^2Q(x_2) \tag{4.10}$$

In order to show the qualitative difference of conditions (4.9) and (4.10) notice that the locally Lipschitz function $g(x):=x_2x_3\left(|x_2|-\gamma(|x_3|)\right)\left(\gamma(|x_3|)-\lambda|x_2|\right)$ satisfies condition (4.10) for arbitrary selection of $\lambda\in(0,1)$ and $\gamma\in K_\infty\cap C^1(\Re_+)$ with $\frac{d\gamma}{ds}(0)>0$. On the other hand, the function $g(x):=x_2x_3\left(|x_2|-\gamma(|x_3|)\right)\left(\gamma(|x_3|)-\lambda|x_2|\right)$ does not satisfy condition (4.9) for any choice of $p,q>0$. To see this, notice that the equation $qx_3=-px_2g(x)$ holds for $x_2\neq 0$ and $x_3=\gamma^{-1}\left((\lambda+1)2^{-1}|x_2|+2^{-1}\sqrt{(\lambda-1)^2x_2^2+4qp^{-1}x_2^{-2}}\right)$. However, the inequality $x_1^2-x_1x_2+px_2^2>qx_1^2x_3$ cannot be satisfied for arbitrary $x_1=1$, $x_2\neq 0$ and $x_3=\gamma^{-1}\left((\lambda+1)2^{-1}|x_2|+2^{-1}\sqrt{(\lambda-1)^2x_2^2+4qp^{-1}x_2^{-2}}\right)$ (notice that if $x_2\to 0$ then $x_3=\gamma^{-1}\left((\lambda+1)2^{-1}|x_2|+2^{-1}\sqrt{(\lambda-1)^2x_2^2+4qp^{-1}x_2^{-2}}\right)\to+\infty$).

In order to obtain condition (4.10), we apply Corollary 4.1 with $\eta(x)\equiv -1$, $\varepsilon=-1$, a function $Q\in C^0(\Re;\Re_+)$ satisfying $0<Q(x)\leq\sigma$ for all $x\neq 0$, where $\sigma\in(0,1)$. We first notice that implications (4.3) and the small-gain conditions (3.10) hold for arbitrary $\gamma\in K_\infty$, $\lambda\in(0,1)$ and the selections:

$$\tilde{\gamma}_{1,2}(s):=\frac{s}{1-\sigma}\quad,\quad \tilde{\gamma}_{1,3}(s):=0 \tag{4.11}$$

$$\tilde{\gamma}_{2,1}(s):=\lambda(1-\sigma)s\quad,\quad \tilde{\gamma}_{2,3}(s):=\gamma(s) \tag{4.12}$$

$$\tilde{\gamma}_{3,1}(s):=\gamma^{-1}\left(\lambda(1-\sigma)s\right)\quad,\quad \tilde{\gamma}_{3,2}(s):=\gamma^{-1}\left(\lambda s\right) \tag{4.13}$$

Since $\gamma\in K_\infty\cap C^1(\Re_+)$ with $\frac{d\gamma}{ds}(0)>0$ there exist constants $0<\gamma_1\leq\gamma_2$ such that

$$\gamma_1 s\leq\gamma(s)\leq\gamma_2 s\quad,\quad\forall s\in[0,r] \tag{4.14}$$

for $r>0$ sufficiently small. We next notice that there exist $p,r>0$ such that the vector $k=(0,0,-p)\in\Re^3$ achieves $-x_2^2-px_2x_3g(x)\leq-\sigma x_2^2$ for all $x\in\Re^3$ with $|x|\leq r$, $\max_{s=1,3}\tilde{\gamma}_{2,s}(|x_s|)\leq|x_2|$ and $x_3x_1^2-px_3^2\leq-\sigma x_3^2$ for all $x\in\Re^3$ with $|x|\leq r$, $\max_{s=1,2}\tilde{\gamma}_{3,s}(|x_s|)\leq|x_3|$. More specifically, the constant $p>0$ must satisfy

$$\frac{r\gamma_2}{\lambda(1-\sigma)}+\sigma<p<\frac{\gamma_1(1-\sigma)}{Lr+|g(0)|} \tag{4.15}$$

where $L\geq 0$ is the Lipschitz constant that satisfies $|g(x)-g(0)|\leq L|x|$, for all $x\in\Re^3$ with $|x|\leq r$. Since $\frac{d\gamma}{ds}(0)>\frac{\sigma|g(0)|}{1-\sigma}$, it follows from (4.14) that a selection of $p>0$ according to (4.15) is possible provided that $r>0$ is selected to be sufficiently small.

Finally, we check implication (4.2). Indeed, using definitions (4.11), (4.12), (4.13), we conclude that implication (4.2) is equivalent to implication (4.10).

The reader should notice that it might be possible to obtain a single CLF for system (4.8) of the form $V(x)=\varphi_1(x_1)+\varphi_2(x_2)+\varphi_3(x_3)$, where $\varphi_i\in C^1(\Re;\Re_+)$ ($i=1,2,3$) are positive definite, radially unbounded functions under condition (4.10) (using the methodologies for the construction of Lyapunov functions in [9] and [17]). However, it should be noted how easily condition (4.10) was obtained from Corollary 4.1. ◁

However, it should be noted that the use of functions $\eta\in C^1(\Re^n;\Re)$, $W\in C^1(\Re^n;[1,+\infty))$, $\delta\in C^0(\Re_+;(0,+\infty))$ and $K\in C^0(\Re_+;[1,+\infty))$ can give much less demanding conditions for the existence of a smooth globally stabilizing feedback. The following example illustrates this point.

**Example 4.4:** Again, we study problem (Q) for system (4.8). Here, we assume that $g:\Re^3 \to \Re$ is independent of $x_3$, i.e., $g(x)=g(x_1,x_2)$. In this example, we will show that system (4.8) can be stabilized globally by smooth feedback provided that there exist constants $\lambda \in (0,1)$, $R>0$, $\sigma \in (0,1)$, functions $\gamma \in K_\infty \cap C^1(\Re_+)$ with $\frac{d\gamma}{ds}(0) > \frac{\sigma |g(0)|}{1-\sigma}$, $q \in C^0(\Re;\Re_+)$ with $0<q(x)\le \sigma$ for all $x \neq 0$, such that the locally Lipschitz function $g:\Re^3 \to \Re$ satisfies the following implication:

$$\lambda(1-\sigma)|x_1| \le \gamma(|x_3|),\ \lambda|x_2| \le \gamma(|x_3|) \le |x_2|,\ x_1^2+x_2^2 \le R,\ x_2x_3g(x)<0 \Rightarrow x_2^2(1-q(x_2))+x_2x_1^2g(x) \ge |x_2x_3g(x)|q(x_3) \tag{4.16}$$

Implication (4.16) is less demanding than implication (4.10) since the inequality $x_2^2(1-Q(x_2))+x_2x_1^2g(x) \ge |x_2x_3g(x)|Q(x_3)$ is assumed to hold only for points that belong to the set $S_1 := \{x\in\Re^3 : \lambda(1-\sigma)|x_1| \le \gamma(|x_3|),\ \lambda|x_2| \le \gamma(|x_3|) \le |x_2|, x_2x_3g(x)<0, x_1^2+x_2^2 \le R\}$. Notice that implication (4.10) requires that the inequality $x_2^2(1-Q(x_2))+x_2x_1^2g(x) \ge |x_2x_3g(x)|Q(x_3)$ holds for points that belong to the set $S_2 := \{x\in\Re^3 : \lambda(1-\sigma)|x_1| \le \gamma(|x_3|),\ \lambda|x_2| \le \gamma(|x_3|) \le |x_2|, x_2x_3g(x)<0\}$ and $S_1 \subseteq S_2$. For example, any function $g:\Re^3\to\Re$ with $g(x)x_2 \ge 0$ for all $x\in\Re^3$ in the compact set $S := \{x\in\Re^3 : \lambda(1-\sigma)|x_1| \le \gamma(|x_3|),\ \lambda|x_2| \le \gamma(|x_3|) \le |x_2|, x_1^2+x_2^2 \le R\}$ satisfies implication (4.16) for appropriate $q \in C^0(\Re;\Re_+)$ but does not necessarily satisfies implication (4.10).

In order to obtain condition (4.16), we apply Corollary 4.1 with $\eta(x) := -a + (x_1^2+x_2^2)/2$, $W(x) := 1+(x_1^2+x_2^2+x_3^2)/2$, $\bar{k}(x) \equiv 0$, $\delta(\eta) \equiv c$, $K(\eta) := 2(\eta+a)+c^{-1}$, $a \ge c > 0$, $\varepsilon > 0$ sufficiently small constants and an appropriate function $Q \in C^0(\Re;\Re_+)$ satisfying $0<Q(x)\le\sigma$ for all $x\neq 0$, where $\sigma\in(0,1)$. By virtue of (4.16), we notice that implications (4.2), (4.3) and the small-gain conditions (3.10) hold for the selections given by (4.11), (4.12) and (4.13), provided that:

$$2(a+\varepsilon) \le R \text{ and } Q(y) \le q(y), \text{ for all } y \in \Re \tag{4.17}$$

Moreover, we notice (exactly as in Example 4.3) that there exist $p,r>0$ such that the vector $k=(0,0,-p)\in\Re^3$ achieves $-x_2^2 - p x_2 x_3 g(x) \le -\sigma x_2^2$ for all $x\in\Re^3$ with $|x|\le r$ with $\max_{s=1,3}\tilde{\gamma}_{2,s}(|x_s|) \le |x_2|$ and $x_3x_1^2 - px_3^2 \le -\sigma x_3^2$ for all $x\in\Re^3$ with $|x|\le r$ with $\max_{s=1,2}\tilde{\gamma}_{3,s}(|x_s|) \le |x_3|$. Furthermore, we obtain for all $x\in\Re^3$ satisfying $\eta(x)\ge 0$:

$$\sum_{i=1}^{3}\frac{\partial W}{\partial x_i}(x)\left(f_i(x)+g_i(x)\bar{k}(x)\right) = -x_1^2+x_1x_2-x_2^2+x_3x_1^2 \le x_3x_1^2 \le 2(\eta(x)+a)W(x) \le K(\eta(x))W(x) \tag{4.18}$$

and

$$\sum_{i=1}^{3}\frac{\partial \eta}{\partial x_i}(x)\left(f_i(x)+g_i(x)\bar{k}(x)\right) = -x_1^2+x_1x_2-x_2^2 \le -\frac{1}{2}x_1^2-\frac{1}{2}x_2^2 = -(\eta(x)+a) \le -a \le -c = -\delta(\eta(x)) \tag{4.19}$$

Finally, we check implications (4.4), (4.5). Implications (4.4) are equivalent to the following implication:

$$\text{"If } \lambda(1-\sigma)|x_1| \le \gamma(|x_3|),\ \lambda|x_2| \le \gamma(|x_3|),\ 2a \le x_1^2+x_2^2 \le 2(a+\varepsilon) \text{ and } x_3x_2g(x)<0 \text{ then } -x_1^2+x_1x_2-x_2^2-x_2g(x)\left(x_1^2+x_3Q(x_3)\right) \le -\delta(\eta(x))\text{"} \tag{4.20}$$

and implications (4.5) are equivalent to the following implications:

$$\text{"If } \lambda(1-\sigma)|x_1| \le |x_2|,\ \gamma(|x_3|) \le |x_2|,\ 2a \le x_1^2+x_2^2 \le 2(a+\varepsilon) \text{ and } x_2g(x)(x_2g(x)+x_3)<0 \text{ then } -x_1^2+x_1x_2-x_2^2+x_3x_1^2-x_2(x_2g(x)+x_3)\frac{Q(x_2)-1}{g(x)} \le K(\eta(x))W(x)\text{"} \tag{4.21}$$

"If $\lambda(1-\sigma)|x_1| \le \gamma(|x_3|)$, $\lambda|x_2| \le \gamma(|x_3|)$, $2a \le x_1^2 + x_2^2 \le 2(a+\varepsilon)$ and $x_3(x_2 g(x) + x_3) < 0$ then $-x_1^2 + x_1 x_2 - x_2^2 - x_2 g(x) x_1^2 - (x_2 g(x) + x_3) x_3 Q(x_3) \le K(\eta(x)) W(x)$" (4.22)

Since $g : \Re^3 \to \Re$ is independent of $x_3$, i.e., $g(x) = g(x_1, x_2)$, it follows from (4.17), (4.20)-(4.22) that if $a, \varepsilon, c > 0$ are sufficiently small constants then all requirements of Corollary 4.1 hold with $\delta(\eta) \equiv c$, $K(\eta) := 2(\eta + a) + c^{-1}$ and appropriate function $Q \in C^0(\Re; \Re_+)$. $\triangleleft$

## 5. Stabilization of Reaction Networks

Reaction networks taking place in Continuous Stirred Tank Reactors are described by ordinary differential equations of the form:

$$\dot{c} = D(c_f - c) + S\upsilon(c) \tag{5.1}$$

where $c \in \mathrm{int}(\Re_+^n)$ is the vector of concentrations of all ($n$) species, $\upsilon : \Re_+^n \to \Re_+^m$ is the continuously differentiable vector field of the ($m$) reaction rates, $S \in \Re^{n \times m}$ is the so-called stoichiometric matrix, $c_f \in \Re_+^n$ is the constant vector of the concentrations of the species at the reactor inlet (feed) and $D \in [0, D_{\max}]$, where $D_{\max} > 0$, is the so-called dilution rate and is the ratio of the volumetric feed rate (equal to the volumetric outlet rate) over the volume of the tank reactor. The dilution rate is used as input in many cases for the achievement of certain control objectives.

In general, the vector field $\upsilon : \Re_+^n \to \Re_+^m$ satisfies the following condition:

$$\text{"If } c_i = 0 \text{ and } S_{i,j} < 0 \text{ then } \upsilon_j(c) = 0\text{"} \tag{5.2}$$

which expresses the (logical) requirement that a reaction cannot occur if one of its reactants is absent.

The equilibrium points of (5.1) for $D = D^* \in (0, D_{\max})$ satisfy the following equation:

$$c_f = c^* - (D^*)^{-1} S\upsilon(c^*) \tag{5.3}$$

Notice that without loss of generality we may assume that $D^* = 1 \in (0, D_{\max})$ (by applying appropriate time-scaling).

This section is devoted to the global stabilization problem of one of the equilibrium points of the reactor. More specifically, we will study the reaction network (5.1) under the following hypotheses:

**(R1)** There exist $N$ pairs of vectors $q_l \in \Re_+^m$, $p_l \in \Re^n$ ($l = 1,...,N$) such that:

$$S' p_l = q_l \quad l = 1,...,N \tag{5.4}$$

**(R2)** There exist constants $b, R > 0$ such that

$$\max_{i=1,...,n}(c_i) \le b + R\sum_{l=1}^{N} \max(p_l' c_f - p_l' c, 0), \text{ for all } c \in \Re_+^n \tag{5.5}$$

**(R3)** There exists $g \in K_\infty$ such that for all $j \in \{1,...,m\}$, $i \in \{1,...,n\}$ with $S_{i,j} < 0$ it holds that

$$0 \le \upsilon_j(c) \le g\left(\max_{j=1,...,n}(c_j)\right) c_i, \text{ for all } c \in \Re_+^n \tag{5.6}$$

The reader should notice that hypothesis (R3) is a direct consequence of the requirement (5.2) and the fact that $\upsilon : \Re_+^n \to \Re_+^m$ is a continuously differentiable vector field. Hypotheses (R1) and (R2) usually hold for all reaction networks. More specifically, the conservation of total mass requires the existence $p \in \Re^n$ such that $S'p = 0$, provided that all chemical species of the reaction mixture are accounted in the description of the network. On the other hand, hypotheses (R1) and (R2) can allow the case where some of the reaction products are not accounted (because they are inert). Hypotheses (R1), (R2) and (R3) usually hold for biochemical networks as well. All chemostat models (see [23]) satisfy hypotheses (R1), (R2) and (R3).

The control problem of the global stabilization of one of the equilibrium points of (5.1) is a meaningful problem because many times there are multiple equilibrium points, indicating absence of global stability. The following example illustrates that the phenomenon of multiple equilibrium points can happen even for very simple reaction networks.

**Example 5.1:** Consider the simple reaction network $1 \to 2$ taking place in a CSTR:

$$\dot{c} = D\left(c_f - c\right) + \begin{bmatrix} -1 \\ 1 \end{bmatrix} k c_1 c_2^2 \tag{5.7}$$

where $c = (c_1, c_2)' \in \mathrm{int}(\Re_+^2)$, $k > 0$ is a constant and $c_f = (c_{1,f}, c_{2,f})' \in \Re_+^2$ with $c_{1,f} > 0$. It is clear that the simple reaction network (5.7) takes the form (5.1) with $n = 2$, $m = 1$, $\upsilon(c) = k c_1 c_2^2$ and $S' = [-1, 1]$. The equilibrium points of (5.7) satisfy the equations:

$$c_{1,f} = c_1^* + k c_1^* \left(c_2^*\right)^2 \quad , \quad c_{2,f} = c_2^* - k c_1^* \left(c_2^*\right)^2 \tag{5.8}$$

The above system of equations has a unique solution if $k\left(c_{1,f} + c_{2,f}\right) < 3$. On the other hand, if $k\left(c_{1,f} + c_{2,f}\right) > 3$ then the above system of equations admits three different solutions. All solutions of (5.8) satisfy $c_{1,f} > c_1^* > 0$, $0 \le c_{2,f} < c_2^* < c_{1,f} + c_{2,f}$ and $c_1^* + c_2^* = c_{1,f} + c_{2,f}$. It is clear that the global stabilization problem for one of the equilibrium points of (5.7) is particularly meaningful for the case $k\left(c_{1,f} + c_{2,f}\right) > 3$.

The reader should notice that hypothesis (R1) holds with $p_1 = \begin{bmatrix} -1 \\ -1 \end{bmatrix}$ and $q_1 = [0] \in \Re$. Other pairs of vectors can be found (e.g. $p_2 = \begin{bmatrix} 0 \\ 1 \end{bmatrix}$ and $q_2 = [1] \in \Re$). Inequality (5.5) with $R = 1$ and $b := c_{1,f} + c_{2,f}$ is a direct consequence of the following inequality:

$$\max_{i=1,2}\left(c_i\right) \le c_1 + c_2 \le c_{1,f} + c_{2,f} + \max\left(p_1' c_f - p_1' c, 0\right), \text{ for all } c \in \Re_+^2 \tag{5.9}$$

Finally, hypothesis (R3) holds with $g(s) := k s^2$. This example will be continued. ◁

Without loss of generality, if there exists an equilibrium point $c^* \in \mathrm{int}(\Re_+^n)$ satisfying (5.3) with $D^* = 1 \in (0, D_{\max})$ then we may assume that $c^* = \mathbf{1}_n$.

The following theorem provides sufficient conditions for the existence of a smooth stabilizing feedback.

**Theorem 5.2:** *Consider system (5.1) under hypotheses (R1), (R2) and (R3). Assume that* $c^* = \mathbf{1}_n \in \mathrm{int}(\Re_+^n)$ *satisfies (5.3) with* $D^* = 1$. *Moreover, suppose that there exist functions* $\tilde{\gamma}_{i,j} \in N_1$, $i, j = 1, \ldots, n$, *with* $\gamma_{i,i}(s) \equiv 0$ *for* $i = 1, \ldots, n$ *which satisfy the small-gain conditions (3.10), a function* $\tilde{Q} \in C^0((0, +\infty); \Re_+)$ *with* $\tilde{Q}(x) > 0$ *for all* $x \in (0,1) \cup (1, +\infty)$ *and constants* $\varepsilon, \omega > 0$ *such that the following implications hold for all* $i, j \in \{1, \ldots, n\}$:

"If $\max\limits_{s=1,\ldots,n} \tilde{\gamma}_{i,s}\left(\left|\ln(c_s)\right|\right) \le \left|\ln(c_i)\right|$, $\max\limits_{s=1,\ldots,n} \tilde{\gamma}_{j,s}\left(\left|\ln(c_s)\right|\right) \le \left|\ln(c_j)\right|$, $\sum\limits_{l=1}^{N}\left(\max\left(p_l' c_f - p_l' c, 0\right)\right)^2 \le 2\varepsilon$

and one of the statements (D1), (D2), (D3) and (D4) holds

$$\text{then } \frac{\left(c_{j,f} - c_j\right)\sum\limits_{l=1}^{m} S_{i,l}\upsilon(c) - \left(c_{i,f} - c_i\right)\sum\limits_{l=1}^{m} S_{j,l}\upsilon(c)}{c_i \ln\left(c_i\right)\tilde{Q}(c_i)\left(c_{j,f} - c_j\right) - c_j \ln\left(c_j\right)\tilde{Q}(c_j)\left(c_{i,f} - c_i\right)} \le -1\text{"} \tag{5.10}$$

"If $c_i = c_{i,f} > 0$, $\sum_{l=1}^{N} \left( \max\left( p_l' c_f - p_l' c, 0 \right) \right)^2 \le 2\varepsilon$ and $\max_{s=1,\ldots,n} \tilde{\gamma}_{i,s}\left( \left| \ln(c_s) \right| \right) \le \left| \ln(c_i) \right|$

$$\text{then } \ln\left(c_{i,f}\right) \sum_{l=1}^{m} S_{i,l}\upsilon_l(c) + c_{i,f}\left(\ln\left(c_{i,f}\right)\right)^2 \tilde{Q}(c_{i,f}) \le 0 \text{"} \tag{5.11}$$

"If $\max_{s=1,\ldots,n} \tilde{\gamma}_{i,s}\left( \left| \ln(c_s) \right| \right) \le \left| \ln(c_i) \right|$, $\varepsilon \le \sum_{l=1}^{N} \left( \max\left( p_l' c_f - p_l' c, 0 \right) \right)^2 \le 2\varepsilon$ and $\min\left(c_{i,f}, 1\right) < c_i < \max\left(c_{i,f}, 1\right)$

$$\text{then } \omega \le -\frac{2\varepsilon}{c_{i,f} - c_i} \left( c_i \ln(c_i) \tilde{Q}(c_i) + \sum_{j=1}^{m} S_{i,j}\upsilon_j(c) \right) \text{"} \tag{5.12}$$

"If $\max_{s=1,\ldots,n} \tilde{\gamma}_{i,s}\left( \left| \ln(c_s) \right| \right) \le \left| \ln(c_i) \right|$, $\sum_{l=1}^{N} \left( \max\left( p_l' c_f - p_l' c, 0 \right) \right)^2 \le 2\varepsilon$ and $\min\left(c_{i,f}, 1\right) < c_i < \max\left(c_{i,f}, 1\right)$

$$\text{then } \ln(c_i) \sum_{l=1}^{m} S_{i,l}\upsilon_l(c) + c_i \left(\ln(c_i)\right)^2 \tilde{Q}(c_i) < 0 \text{"} \tag{5.13}$$

"If $\max_{s=1,\ldots,n} \tilde{\gamma}_{i,s}\left( \left| \ln(c_s) \right| \right) \le \left| \ln(c_i) \right|$, $\sum_{l=1}^{N} \left( \max\left( p_l' c_f - p_l' c, 0 \right) \right)^2 \le 2\varepsilon$ and $c_i > \max\left(c_{i,f}, 1\right)$ or $\min\left(c_{i,f}, 1\right) > c_i$

$$\text{then } \ln(c_i) \sum_{l=1}^{m} S_{i,l}\upsilon_l(c) + c_i \left(\ln(c_i)\right)^2 \tilde{Q}(c_i) + D_{\max}\left(c_{i,f} - c_i\right) \ln(c_i) < 0 \text{"} \tag{5.14}$$

*where the logical statements (D1), (D2), (D3) and (D4) are expressed by*

**(D1)** $\min\left(c_{i,f}, 1\right) < c_i < \max\left(c_{i,f}, 1\right)$ and $c_j > \max\left(c_{j,f}, 1\right)$
**(D2)** $\min\left(c_{i,f}, 1\right) < c_i < \max\left(c_{i,f}, 1\right)$ and $\min\left(c_{j,f}, 1\right) > c_j$
**(D3)** $c_i > \max\left(c_{i,f}, 1\right)$ and $\min\left(c_{j,f}, 1\right) < c_j < \max\left(c_{j,f}, 1\right)$
**(D4)** $\min\left(c_{i,f}, 1\right) > c_i$ and $\min\left(c_{j,f}, 1\right) < c_j < \max\left(c_{j,f}, 1\right)$

*Moreover, suppose that there exists a vector* $k = (k_1, \ldots, k_n)' \in \Re^n$ *such that* $\left( 1 + \sum_{j=1}^{n} k_j \ln\left(c_j\right) \right) \ln(c_i)\left(c_{i,f} - c_i\right) + \ln(c_i) \sum_{j=1}^{m} S_{i,j}\upsilon_j(c) \le -\left(\ln(c_i)\right)^2 c_i \tilde{Q}(c_i)$ *for all* $c \in \mathrm{int}(\Re_+^n)$ *in a neighborhood of* $\mathbf{1}_n$ *with* $\max_{s=1,\ldots,n} \tilde{\gamma}_{i,s}\left( \left| \ln(c_s) \right| \right) \le \left| \ln(c_i) \right|$ *(* $i = 1, \ldots, n$ *).*

*Then there exists* $\tilde{k} \in C^\infty(\mathrm{int}(\Re_+^n); [0, D_{\max}])$, *with* $\tilde{k}(\mathbf{1}_n) = 1$, *such that* $\mathbf{1}_n \in \mathrm{int}(\Re_+^n)$ *is UGAS for the closed-loop system (5.1) with* $D = \tilde{k}(c)$.

**Proof:** The proof utilizes Corollary 4.2 for the system obtained by the following change of coordinates

$$c = \exp(x) \tag{5.15}$$

and the input transformation

$$D = 1 + u \tag{5.16}$$

namely, the system:

$$\dot{x}_i = (1+u)\left(c_{i,f} e^{-x_i} - 1\right) + e^{-x_i} \sum_{j=1}^{m} S_{i,j} \upsilon_j(e^x) \,,\; i = 1,...,n \tag{5.17}$$

It follows from (5.16) that $u \in U := [-1, D_{\max} - 1]$ and consequently hypothesis (P3) holds for system (5.17). Moreover, $D_{\max} - 1 > 0$. Define

$$P(i) := \left\{ j \in \{1,...,m\} : S_{i,j} < 0 \right\} \text{ and } \sigma := \max_{i=1,...,n} \left( \sum_{j \in P(i)} \left| S_{i,j} \right| \right) \tag{5.18}$$

We will apply Corollary 4.2 with $f_i(x) = c_{i,f} e^{-x_i} - 1 + e^{-x_i} \sum_{j=1}^{m} S_{i,j} \upsilon_j(e^x)$, $g_i(x) = c_{i,f} e^{-x_i} - 1$ for $i = 1,...,n$ and

$$\eta(x) := \sum_{l=1}^{N} \left( \max\left( p_l' c_f - p_l' e^x , 0\right)\right)^2 - \varepsilon = \sum_{l=1}^{N} \left( \max\left( p_l' c_f - p_l' c , 0\right)\right)^2 - \varepsilon \tag{5.19}$$

$$W(x) := 1 + \varepsilon + \eta(x) + \sum_{i=1}^{n} e^{-x_i} \tag{5.20}$$

$$\bar{k}(x) \equiv 0 \,,\; \delta(\eta) \equiv \min(\varepsilon, \omega) \,,\; Q(x) := \tilde{Q}\left(e^x\right) \tag{5.21}$$

$$K(\eta) := D_{\max} + \sigma g\left( b + (N+\varepsilon)\frac{R}{2} + \frac{R}{2}\eta \right) \tag{5.22}$$

All functions defined above satisfy the requirements imposed by Corollary 4.2. More specifically, using (5.3) and (5.4) we get $p_l' c_f - p_l' \mathbf{1}_n = -q_l' \upsilon(\mathbf{1}_n)$ for $l = 1,...,N$. Since $q_l \in \Re_+^m$ and $\upsilon(\mathbf{1}_n) \in \Re_+^m$, it follows that $q_l' \upsilon(\mathbf{1}_n) \geq 0$ and $p_l' c_f \leq p_l' \mathbf{1}_n$. Consequently, definition (5.19) implies $\eta(0) = -\varepsilon < 0$. Next notice that definitions (5.19), (5.20) imply that $W(x) \geq 1$, $\max_{i=1,...,n} \left(-x_i\right) \leq \ln\left(W(x)\right)$ for all $x \in \Re^n$. Inequality (5.5) implies $\max_{i=1,...,n} \left(e^{x_i}\right) \leq b + N\frac{R}{2} + \frac{R}{2} \sum_{l=1}^{N} \left( \max\left( p_l' c_f - p_l' e^x , 0\right)\right)^2$ for all $x \in \Re^n$, which combined with definition (5.19) gives:

$$\max_{i=1,...,n} \left(e^{x_i}\right) \leq b + N\frac{R}{2} + \frac{R}{2}(\eta(x) + \varepsilon) \tag{5.23}$$

The above inequalities in conjunction with definition (5.20), allow us to conclude that the following inequalities hold for all $x \in \Re^n$:

$$\max_{i=1,...,n} \left(-x_i\right) \leq \ln\left(W(x)\right),\; \max_{i=1,...,n} \left(x_i\right) \leq \ln\left( b + N\frac{R}{2} + \frac{R}{2} W(x) \right) \tag{5.24}$$

Inequality (5.24) shows that $W \in C^1(\Re^n; [1,+\infty))$ as defined by (5.20) is radially unbounded. Using (5.4) and definitions (5.19), (5.20) we obtain by differentiating $\eta, W$ for all $x \in \Re^n$ and $u \in [-1, D_{\max} - 1]$:

$$\sum_{i=1}^{n} \frac{\partial \eta}{\partial x_i}(x)\left(f_i(x) + g_i(x)u\right) = -2(1+u)(\varepsilon + \eta(x)) - 2\sum_{l=1}^{N} \left( \max\left( p_l' c_f - p_l' e^x , 0\right)\right) q_l' \upsilon(e^x) \tag{5.25}$$

$$\sum_{i=1}^{n}\frac{\partial W}{\partial x_i}(x)\big(f_i(x)+g_i(x)u\big)=-(1+u)\left(2\big(\varepsilon+\eta(x)\big)+\sum_{i=1}^{n}e^{-2x_i}\left(c_{i,f}-e^{x_i}\right)\right)$$
$$-2\sum_{l=1}^{N}\Big(\max\big(p_l'c_f-p_l'e^x,0\big)\Big)q_l'\upsilon(e^x)-\sum_{i=1}^{n}e^{-2x_i}\sum_{j=1}^{m}S_{i,j}\upsilon_j(e^x) \tag{5.26}$$

Since $q_l\in\Re_+^m$ and $\upsilon(e^x)\in\Re_+^m$, it follows that $\Big(\max\big(p_l'c_f-p_l'e^x,0\big)\Big)q_l'\upsilon(e^x)\geq 0$ for $l=1,...,N$. Therefore, we obtain from (5.25) and (5.26) for all $x\in\Re^n$ and $u\in[-1,D_{\max}-1]$:

$$\sum_{i=1}^{n}\frac{\partial\eta}{\partial x_i}(x)\big(f_i(x)+g_i(x)u\big)\leq -2(1+u)(\varepsilon+\eta(x)) \tag{5.27}$$

$$\sum_{i=1}^{n}\frac{\partial W}{\partial x_i}(x)\big(f_i(x)+g_i(x)u\big)\leq D_{\max}\sum_{i=1}^{n}e^{-x_i}-\sum_{i=1}^{n}e^{-2x_i}\sum_{j=1}^{m}S_{i,j}\upsilon_j(e^x) \tag{5.28}$$

Using (5.6) and (5.18) we get $-\sum_{i=1}^{n}e^{-2x_i}\sum_{j=1}^{m}S_{i,j}\upsilon_j(e^x)\leq\sum_{i=1}^{n}e^{-x_i}\sum_{j\in P(i)}\left|S_{i,j}\right|g\left(\max_{i=1,...,n}\left(e^{x_i}\right)\right)$ for all $x\in\Re^n$. Hence, we obtain from (5.18), (5.20), (5.22), (5.23) and (5.28) for all $x\in\Re^n$ and $u\in[-1,D_{\max}-1]$:

$$\sum_{i=1}^{n}\frac{\partial W}{\partial x_i}(x)\big(f_i(x)+g_i(x)u\big)\leq\left(D_{\max}+\sigma g\left(b+(N+\varepsilon)\frac{R}{2}+\frac{R}{2}\eta(x)\right)\right)\sum_{i=1}^{n}e^{-x_i}$$
$$\leq\left(D_{\max}+\sigma g\left(b+(N+\varepsilon)\frac{R}{2}+\frac{R}{2}\eta(x)\right)\right)W(x)\leq K(\eta(x))W(x) \tag{5.29}$$

Using (5.21), (5.27) and (5.29) we conclude that the function $\bar{k}\in C^0\left(\Re^n;U\right)$ defined by (5.21) satisfies $\sum_{i=1}^{n}\frac{\partial W}{\partial x_i}(x)\big(f_i(x)+g_i(x)\bar{k}(x)\big)\leq K\big(\eta(x)\big)W(x)$ and $\sum_{i=1}^{n}\frac{\partial\eta}{\partial x_i}(x)\big(f_i(x)+g_i(x)\bar{k}(x)\big)\leq-\delta\big(\eta(x)\big)$ for all $x\in\Re^n$ with $\eta(x)\geq 0$.

We next notice that by virtue of the change of coordinates (5.15), definition (5.19) and the definition of $Q$ in (5.21), it follows that implications (4.2), (4.3), (4.6) and (4.7) are equivalent to implications (5.10), (5.11), (5.13) and (5.14), respectively. Moreover, implication (5.12) implies implication (4.4). Indeed, first notice that (5.25) implies $\sum_{i=1}^{n}\frac{\partial\eta}{\partial x_i}(x)g_i(x)=-2(\varepsilon+\eta(x))<0$ for all $x\in\Re^n$ with $\eta(x)\geq 0$ and consequently the condition $x_jg_j(x)\sum_{l=1}^{n}\frac{\partial\eta}{\partial x_l}(x)g_l(x)<0$ implies $x_jg_j(x)>0$, or equivalently, $\min\left(c_{i,f},1\right)<c_i<\max\left(c_{i,f},1\right)$. Implication (4.4) requires $\sum_{l=1}^{n}\frac{\partial\eta}{\partial x_l}(x)\big(f_l(x)+g_l(x)u\big)\leq-\delta\big(\eta(x)\big)$ for $u=-\frac{f_j(x)+x_jQ(x_j)}{g_j(x)}=-1-\frac{e^{-x_j}\sum_{l=1}^{m}S_{j,l}\upsilon_l(e^x)+x_jQ(x_j)}{c_{j,f}e^{-x_j}-1}$ or $1+u=-\frac{\sum_{l=1}^{m}S_{j,l}\upsilon_l(c)+c_j\tilde{Q}(c_j)\ln\left(c_j\right)}{c_{j,f}-c_j}$. Notice that (5.12) implies that $\frac{\omega}{2\varepsilon}\leq-\frac{c_j\ln\left(c_j\right)\tilde{Q}(c_j)+\sum_{l=1}^{m}S_{j,l}\upsilon_l(c)}{c_{j,f}-c_j}$ and inequality (5.27) combined with definition (5.21) gives $\sum_{i=1}^{n}\frac{\partial\eta}{\partial x_i}(x)\big(f_i(x)+g_i(x)u\big)\leq-\omega\leq-\delta\big(\eta(x)\big)$ for all $x\in\Re^n$ with $0\leq\eta(x)\leq\varepsilon$. Consequently, implication (4.4) holds.

We next show that implication (4.5) holds. By virtue of (5.29) it suffices to show that for all $x \in \Re^n$ with $\max_{s=1,\dots,n} \tilde{\gamma}_{j,s}(|x_s|) \le |x_j|$, $0 \le \eta(x) \le \varepsilon$:

- $-\dfrac{f_j(x) + x_j Q(x_j)}{g_j(x)} \le D_{\max} - 1$ for $x_j g_j(x) < 0$ and $\sum_{l=1}^{n} \dfrac{\partial W}{\partial x_l}(x) g_l(x) > 0$,
- $-\dfrac{f_j(x) + x_j Q(x_j)}{g_j(x)} \ge -1$ for $x_j g_j(x) > 0$ and $\sum_{l=1}^{n} \dfrac{\partial W}{\partial x_l}(x) g_l(x) < 0$

The above conditions are implied by the following conditions for all $c \in \mathrm{int}(\Re_+^n)$ with $\max_{s=1,\dots,n} \tilde{\gamma}_{i,s}(|\ln(c_s)|) \le |\ln(c_i)|$, $\sum_{l=1}^{N} \left(\max(p_l' c_f - p_l' c, 0)\right)^2 \le 2\varepsilon$:

- $-\dfrac{\sum_{l=1}^{m} S_{j,l} \upsilon_l(c) + c_j \tilde{Q}(c_j) \ln(c_j)}{c_{j,f} - c_j} \le D_{\max}$ for $c_j > \max(c_{j,f}, 1)$ or $\min(c_{j,f}, 1) > c_j$,
- $0 \le -\dfrac{\sum_{l=1}^{m} S_{j,l} \upsilon_l(c) + c_j \tilde{Q}(c_j) \ln(c_j)}{c_{j,f} - c_j}$ for $\min(c_{j,f}, 1) < c_j < \max(c_{j,f}, 1)$

The above conditions are direct consequences of (5.13) and (5.14). Thus we conclude that implication (4.5) holds.

The proof is complete. $\triangleleft$

The following example illustrates how Theorem 5.2 can be applied to reaction networks. Theorem 5.2 guarantees the existence of a smooth globally stabilizing bounded feedback.

**Example 5.1(continued):** We now turn to the global stabilization problem of one of the equilibrium points $c^* \in \mathrm{int}(\Re_+^2)$ of system (5.7) by means of smooth bounded feedback. In order to apply Theorem 5.2 we first apply a coordinate change that "brings" the equilibrium point to $\mathbf{1}_2 \in \mathrm{int}(\Re_+^2)$ (namely, we apply the coordinate change $c_1 \to c_1^* c_1$, $c_2 \to c_2^* c_2$). Then system (5.7) takes the form

$$\begin{aligned} \dot{c} &= D\left( \begin{bmatrix} 1+\theta \\ 1-\mu\theta \end{bmatrix} - c \right) + \begin{bmatrix} -1 \\ \mu \end{bmatrix} \theta c_1 c_2^2 \\ c &= (c_1, c_2)' \in \mathrm{int}(\Re_+^2), D \in [0, D_{\max}] \end{aligned} \tag{5.30}$$

where $\theta > 0$ and $\mu \in (0, \theta^{-1}]$ are constant parameters ($k\left(c_2^*\right)^2 = \theta$, $\mu = c_1^* / c_2^*$). For the above system $\mathbf{1}_2 \in \mathrm{int}(\Re_+^2)$ is the equilibrium point to be globally stabilized. Hypotheses (R1)-(R3) hold with $N = 1$, $q_1 = [0] \in \Re_+^1$, $p_1 = (-\mu, -1)' \in \Re^2$. We will further assume that $\mu < \theta^{-1}$, i.e., $c_{2,f} > 0$ for system (5.7).

The problem that we will study in this example is:

> "How large must $D_{\max} > 0$ be so that $\mathbf{1}_2 \in \mathrm{int}(\Re_+^2)$ can be globally stabilized by a smooth feedback law $D = k(c)$ with $k(\mathbf{1}_2) = 1$ and $k(c) \in [0, D_{\max}]$ for all $c \in \mathrm{int}(\Re_+^2)$?"

In order to solve the above problem we exploit Theorem 5.2. We will show that all conditions of Theorem 5.2 hold provided that

$$D_{\max} > (1+\mu)^2 \text{ and } D_{\max} > \frac{(1+\mu)\left(\theta(1+\mu)^2+1\right)}{\mu\theta} \tag{5.31}$$

We first check conditions (5.10)-(5.14) of Theorem 5.2. Let $\varepsilon > 0$ be an arbitrary constant, $\gamma \in K_\infty$ be an arbitrary function and select $\tilde{\gamma}_{1,2}(s) := \gamma(s)$, $\tilde{\gamma}_{2,1}(s) := \gamma^{-1}(\lambda s)$ where $\lambda \in (0,1)$ is to be selected (notice that the small-gain conditions are automatically satisfied). Implication (5.10) holds provided that the inequality

$$c_1(\ln(c_1))^2 \tilde{Q}(c_1) - c_2(\ln(c_2))^2 \tilde{Q}(c_2) \frac{\ln(c_1)(1+\theta-c_1)}{\ln(c_2)(1-\mu\theta-c_2)} \leq \theta c_1 c_2^2 \ln(c_1) \frac{\mu+1-\mu c_1 - c_2}{1-\mu\theta-c_2}$$

holds for all $c \in \operatorname{int}(\Re_+^2)$ with $\gamma(|\ln(c_2)|) \leq |\ln(c_1)|$, $\lambda|\ln(c_1)| \leq \gamma(|\ln(c_2)|)$, $\mu c_1 + c_2 \leq 1+\mu+\sqrt{2\varepsilon}$ which satisfy one of the logical statements (D1), (D2), (D3) and (D4). It is clear that the above inequality holds provided that the inequalities

$$(\ln(c_1))^2 \tilde{Q}(c_1) \leq \frac{\theta}{2} c_2^2 \ln(c_1) \frac{\mu+1-\mu c_1 - c_2}{1-\mu\theta-c_2}$$

$$(\ln(c_2))^2 \tilde{Q}(c_2) \leq \frac{\theta}{2} c_1 c_2 \ln(c_2) \frac{\mu c_1 + c_2 - \mu - 1}{1+\theta-c_1}$$

hold for all $c \in \operatorname{int}(\Re_+^2)$ with $\gamma(|\ln(c_2)|) \leq |\ln(c_1)|$, $\lambda|\ln(c_1)| \leq \gamma(|\ln(c_2)|)$, $\mu c_1 + c_2 \leq 1+\mu+\sqrt{2\varepsilon}$ which satisfy one of the logical statements (D1), (D2), (D3) and (D4). The above requirements imply that:

- If (D1) holds, i.e., $1 < c_1 < 1+\theta$ and $1 < c_2 \leq \exp\left(\gamma^{-1}(\ln(c_1))\right)$, then we must have $\tilde{Q}(c_1) \leq \frac{c_1 - 1}{2\ln(c_1)}$ and $\tilde{Q}(c_2) \leq \frac{c_2 - 1}{2\ln(c_2)} c_2$.
- If (D2) holds, i.e., $1 < c_1 < 1+\theta$ and $\exp\left(-\gamma^{-1}(\ln(c_1))\right) \leq c_2 < 1-\mu\theta$, then we must have $\tilde{Q}(c_1) \leq \frac{\theta}{2\ln(c_1)} \exp\left(-2\gamma^{-1}(\ln(c_1))\right)$ and $\tilde{Q}(c_2) \leq \frac{1-c_2}{2|\ln(c_2)|} c_2$.
- If (D3) holds, i.e., $c_1 > 1+\theta$ and $1-\mu\theta < c_2 < 1$, then we must have $\tilde{Q}(c_1) \leq \frac{c_1 - 1}{2\ln(c_1)}(1-\mu\theta)^2$ and $\tilde{Q}(c_2) \leq (1+\theta)\frac{1-c_2}{2|\ln(c_2)|} c_2$.
- If (D4) holds, i.e., $1-\mu\theta < c_2 < 1$ and $\exp\left(-\lambda^{-1}\gamma(|\ln(c_2)|)\right) \leq c_1 < 1$, then we must have $\tilde{Q}(c_1) \leq \frac{1-c_1}{2|\ln(c_1)|}(1-\mu\theta)^2$ and $\tilde{Q}(c_2) \leq \frac{1-c_2}{2|\ln(c_2)|} c_2 \exp\left(-\lambda^{-1}\gamma(|\ln(c_2)|)\right)$.

Consequently, implication (5.10) automatically holds, if $\tilde{Q} \in C^0((0,+\infty);\Re_+)$ is selected to be

$$\tilde{Q}(1) := A\frac{(1-\mu\theta)^2}{2} \tag{5.32}$$

$$\tilde{Q}(c) := A\frac{\min(1,|1-c|)}{2|\ln(c)|}\min\left(c,(1-\mu\theta)^2\right)\exp\left(-\max\left(\lambda^{-1}\gamma(|\ln(c)|), 2\gamma^{-1}(|\ln(c)|)\right)\right), \text{ for } c \in (0,1)\cup(1,+\infty) \tag{5.33}$$

where $A \in (0,1]$ is a constant (yet to be selected).

Conditions (5.11) with $i=1,2$ give the inequalities $\tilde{Q}(1+\theta) \le \frac{\theta}{\ln(1+\theta)} \exp\left(-2\gamma^{-1}\left(\ln(1+\theta)\right)\right)$ and $\tilde{Q}(1-\mu\theta) \le \frac{\mu\theta(1-\mu\theta)}{\left|\ln(1-\mu\theta)\right|} \exp\left(-\lambda^{-1}\gamma\left(\left|\ln(1-\mu\theta)\right|\right)\right)$, which hold automatically for $\tilde{Q} \in C^0((0,+\infty);\Re_+)$ defined by (5.32), (5.33).

Conditions (5.12) with $i=1,2$ and $\omega > 0$ hold provided that

$$\begin{aligned} &A \le \frac{1}{2(1+\theta)} \exp\left(-2\gamma^{-1}\left(\ln(1+\theta)\right)\right) \\ &\omega \le \frac{\varepsilon\theta}{1+\theta} \exp\left(-2\gamma^{-1}\left(\ln(1+\theta)\right)\right) \end{aligned} \quad \text{and} \quad \begin{aligned} &A \le \exp\left(-\lambda^{-1}\gamma\left(\left|\ln(1-\mu\theta)\right|\right)\right) \\ &\omega \le \varepsilon(1-\mu\theta)^2 \exp\left(-\lambda^{-1}\gamma\left(\left|\ln(1-\mu\theta)\right|\right)\right) \end{aligned} \tag{5.34}$$

and conditions (5.13) hold automatically for the above selections. Finally, we check condition (5.14). Conditions (5.14) hold provided that

$$D_{\max} > c_2^2 \text{ and } D_{\max} > \frac{\mu\theta c_1 c_2^2 + c_2 \ln(c_2)\tilde{Q}(c_2)}{\mu\theta} \text{ for all } c \in \operatorname{int}(\Re_+^2) \text{ with } \mu c_1 + c_2 \le (1+\mu) + \sqrt{2\varepsilon}$$

Since $\varepsilon > 0$ is arbitrary, we conclude (by virtue of (5.31)) that the above inequalities hold.

The existence of a vector $k = (k_1, k_2)' \in \Re^2$ with

$$\left(1 + \sum_{j=1}^{2} k_j \ln(c_j)\right) \ln(c_i)(c_{i,f} - c_i) + \ln(c_i) \sum_{j=1}^{m} S_{i,j} \upsilon_j(c) \le -\left(\ln(c_i)\right)^2 c_i \tilde{Q}(c_i)$$

for all $c \in \operatorname{int}(\Re_+^2)$ in a neighborhood of $\mathbf{1}_n$ with $\max_{s=1,\dots,n} \tilde{\gamma}_{i,s}\left(\left|\ln(c_s)\right|\right) \le \left|\ln(c_i)\right|$ ($i=1,2$) is guaranteed by the observation that the unbounded smooth feedback law $D = c_2^2$ guarantees the inequalities $D \ln(c_i)(c_{i,f} - c_i) + \ln(c_i) \sum_{j=1}^{m} S_{i,j} \upsilon_j(c) \le -\left(\ln(c_i)\right)^2 c_i \tilde{Q}(c_i)$ for appropriate functions $\tilde{\gamma}_{1,2} \in K_\infty$, $\tilde{\gamma}_{2,1} \in K_\infty$ and for all $c \in \operatorname{int}(\Re_+^2)$ in a neighborhood of $\mathbf{1}_n$ with $\max_{s=1,\dots,n} \tilde{\gamma}_{i,s}\left(\left|\ln(c_s)\right|\right) \le \left|\ln(c_i)\right|$ ($i=1,2$). ⊲

## 6. Concluding Remarks

This paper has showed how recent results on vector Lyapunov functions can be used for smooth globally stabilizing feedback design for nonlinear systems. In particular, Theorem 3.4 provides necessary and sufficient conditions based on vector control Lyapunov functions for the existence of a smooth global stabilizer for affine in the control uncertain nonlinear systems.

The flexibility of vector Lyapunov functions is a feature that can be exploited for feedback design in large scale systems. Corollaries 4.1 and 4.2 provide simple and easily checkable sufficient conditions for the existence of a smooth global stabilizer for nonlinear systems. Corollaries 4.1 and 4.2 are direct applications of Theorem 3.4 and show how vector Lyapunov functions can lead to results which are not easily obtained by the classical single Lyapunov analysis. The obtained results are exploited for the derivation of sufficient conditions which guarantee stabilizability of the equilibrium point of a reaction network taking place in a continuous stirred tank reactor.

Future research may include:

i) the extension of the present results to the multiple-input case,
ii) the development of explicit formulas for the feedback stabilizers which are designed based on VRCLFs,
iii) the development of "adding an integrator"-like results based on VRCLFs, which can allow important modifications to the backstepping methodology.

# Appendix

**Proof of Lemma 3.5:** Define $B^0 := \{ i \in \{1,...,m\} : g_i = 0 \}$, $B^+ := \{ i \in \{1,...,m\} : g_i > 0 \}$ and $B^- := \{ i \in \{1,...,m\} : g_i < 0 \}$. Notice that implication (I) gives:

$$\text{“If } B^0 \neq \emptyset \text{ then } \max_{i \in B^0}(f_i) < 0\text{”} \tag{A.1}$$

and implication (II) gives:

$$\text{“If } B^+ \neq \emptyset \text{ and } B^- \neq \emptyset \text{ then } \max_{i \in B^+}\left(\frac{f_i}{g_i}\right) < \min_{i \in B^-}\left(\frac{f_i}{g_i}\right)\text{”} \tag{A.2}$$

Consider all possible cases:

(a) $B^+ = \emptyset = B^-$. In this case (A.1) implies that we can select arbitrary $u \in \Re$ so that $f_i + g_i u < 0$ for all $i = 1,...,m$.

(b) $B^+ = \emptyset \neq B^-$. In this case, we can select $u > \max_{i \in B^-}\left(-\frac{f_i}{g_i}\right)$. The previous inequality in conjunction with (A.1) implies that $f_i + g_i u < 0$ for all $i = 1,...,m$.

(c) $B^- = \emptyset \neq B^+$. In this case, we can select $u < \min_{i \in B^+}\left(-\frac{f_i}{g_i}\right)$. The previous inequality in conjunction with (A.1) implies that $f_i + g_i u < 0$ for all $i = 1,...,m$.

(d) $B^+ \neq \emptyset \neq B^-$. In this case (A.2) implies $\max_{i \in B^+}\left(\frac{f_i}{g_i}\right) < \min_{i \in B^-}\left(\frac{f_i}{g_i}\right)$. Therefore, the following inequality holds:

$$\min_{i \in B^+}\left(-\frac{f_i}{g_i}\right) = -\max_{i \in B^+}\left(\frac{f_i}{g_i}\right) > -\min_{i \in B^-}\left(\frac{f_i}{g_i}\right) = \max_{i \in B^-}\left(-\frac{f_i}{g_i}\right)$$

In this case, we can select $u \in \Re$ so that $\min_{i \in B^+}\left(-\frac{f_i}{g_i}\right) > u > \max_{i \in B^-}\left(-\frac{f_i}{g_i}\right)$. The previous inequality in conjunction with (A.1) implies that $f_i + g_i u < 0$ for all $i = 1,...,m$.

Notice that the above selections in each case guarantee that there exists $u \in (-a,+\infty)$ so that $f_i + g_i u < 0$ for all $i = 1,...,m$, provided that either $B^+ = \emptyset$ or $B^+ \neq \emptyset$ and $\min_{i \in B^+}\left(-\frac{f_i}{g_i}\right) > -a$. Implication (III) guarantees that one of the previously mentioned cases holds.

Finally, notice that the above selections in each case guarantee that there exists $u \in (-a,b)$ so that $f_i + g_i u < 0$ for all $i = 1,...,m$, provided that the following implications hold: (i) if $B^+ \neq \emptyset$ then $\min_{i \in B^+}\left(-\frac{f_i}{g_i}\right) > -a$, (ii) if $B^- \neq \emptyset$ then $\max_{i \in B^-}\left(-\frac{f_i}{g_i}\right) < b$. Implications (III) and (IV) guarantee that the previous implications hold.

The converse statements are proved in the same way by distinguishing the above cases. The proof is complete. $\lhd$

**Proof of Lemma 3.6:** The methodology of the proof is to show that every $\gamma_{i,j} \in N_1$ ($i,j = 1,...,k$, $i \neq j$) can be replaced by a function $\tilde{\gamma}_{i,j} \in N_1$ which is positive definite and satisfies $\lim_{s \to +\infty} \tilde{\gamma}_{i,j}(s) = +\infty$ in such a way that the "new" set of functions $\gamma_{p,q} \in N_1$ ($p,q = 1,...,k$) with $\gamma_{i,j} \in N_1$ replaced by $\tilde{\gamma}_{i,j} \in N_1$ satisfies all properties of the VRCLF.

A key observation is that all inequalities (3.3), (3.4), (3.8), (3.9), (3.11), (3.12) and (3.15) hold automatically if $\gamma_{i,j} \in N_1$ is replaced by a function $\tilde{\gamma}_{i,j} \in N_1$ satisfying the following inequality:

$$\gamma_{i,j}(s) \le \tilde{\gamma}_{i,j}(s), \text{ for all } s \ge 0 \tag{A.3}$$

The only thing that remains to be checked is the set of small-gain inequalities (3.10). The inequalities (3.10) which are affected by the replacement of $\gamma_{i,j} \in N_1$ with the function $\tilde{\gamma}_{i,j} \in N_1$ can be expressed by the following inequality:

$$\tilde{\gamma}_{i,j}(a(s)) < s, \text{ for all } s > 0 \tag{A.4}$$

where $a \in N_1$ is defined as follows:

"$a(s)$ for a given $s \ge 0$ is the maximum of $\gamma_{j,i}(s)$ and the maximum of all $(\gamma_{i,z_1} \circ \gamma_{z_1,z_2} \circ \ldots \circ \gamma_{z_{l-1},z_l} \circ \gamma_{z_l,j})(s)$ over the set of all indices $\{z_1,...,z_l\} \subseteq \{1,...,k\} \setminus \{i,j\}$, with $z_p \ne z_q$ if $p \ne q$"

To see this, notice that each inequality (3.10) which is affected by the replacement of $\gamma_{i,j} \in N_1$ with the function $\tilde{\gamma}_{i,j} \in N_1$ is guaranteed by (A.4) and the following fact:

"If $\gamma_1(\gamma_2(s)) < s$ for all $s > 0$ for a given pair of $\gamma_1, \gamma_2 \in N_1$ then $\gamma_2(\gamma_1(s)) < s$ for all $s > 0$."

Notice that since inequalities (3.10) hold with the original $\gamma_{i,j} \in N_1$, it holds that:

$$\gamma_{i,j}(a(s)) < s, \text{ for all } s > 0 \tag{A.5}$$

Define next:

$$\tilde{\gamma}_{i,j}(s) := \max\left(\gamma_{i,j}(s), \frac{1}{2}\tilde{a}^{-1}(s)\right), \text{ for all } s \ge 0 \tag{A.6}$$

where $\tilde{a}(s) := a(s) + s$. Notice that $\tilde{a} \in K_\infty$ and therefore $\tilde{\gamma}_{i,j} \in N_1$ is well defined by definition (A.6). Definition (A.6) and inequality (A.5) guarantee that inequalities (A.3) and (A.4) hold. Moreover, since $\tilde{a}^{-1} \in K_\infty$ it follows that $\tilde{\gamma}_{i,j} \in N_1$ and satisfies $\lim_{s \to +\infty} \tilde{\gamma}_{i,j}(s) = +\infty$. The proof is complete. ◁